\documentclass[10pt]{article}

\usepackage[dvips]{graphicx}
\usepackage{makeidx}
\usepackage{amsmath}
\usepackage{amsfonts}
\usepackage[russian,english]{babel}

\newcommand {\sign}{{\mathrm{sign}}}

\newcommand {\alphavec}{\boldsymbol{\alpha}}

\newcommand {\phivec}{\mbox{\boldmath $\phi$}}

\newcommand {\thetavec}{{\boldsymbol{\theta}}}

\newfont{\pseudocode}{cmtt10}

\newtheorem{thm}{Theorem}
\newtheorem{cor}{Corollary}

\newtheorem{defn}{Definition}

\newtheorem{assume}{Assumption}

\newcommand{\eps}{\varepsilon}

\newcommand{\s}{\mathcal{S}}

\newcommand{\bfx}{\mathbf{x}}
\newcommand{\bff}{\mathbf{f}}

\newcommand{\bfc}{\mathbf{c}}

\newcommand{\bfg}{\mathbf{g}}

\newcommand{\bfs}{\mathbf{s}}

\newcommand{\dpsi}{{\dot\psi}}
\newcommand{\pd}{{\partial}}

\newcommand{\Real}{\mathbb{R}}

\makeindex

\begin{document}


\baselineskip 0.7cm

\title{\bf Parameter estimation and control for a class of systems with nonlinear parametrization}
\author{Ivan Tyukin\thanks{Laboratory for Perceptual Dynamics, RIKEN (Institute for Physical and Chemical Research)
                           Brain Science Institute, 2-1, Hirosawa, Wako-shi, Saitama, 351-0198, Japan, e-mail:
                           \{tyukinivan\}@brain.riken.jp}
                           , Danil
                           Prokhorov\thanks{Ford Research Laboratory, Dearborn, MI,
48121, USA, e-mail: dprokhor@ford.com},  Cees van
Leeuwen\thanks{Laboratory for Perceptual Dynamics, RIKEN
(Institute for Physical and Chemical Research)
                           Brain Science Institute, 2-1, Hirosawa, Wako-shi, Saitama, 351-0198, Japan, e-mail:
                           \{ceesvl\}@brain.riken.jp}
}
\date{\today}
\date{}
\maketitle{}

\begin{abstract}
We propose novel parameter estimation algorithms for a class of
dynamical systems with nonlinear parametrization. The class is
initially restricted to smooth monotonic functions with respect to
a linear functional of the parameters. We show that under this
restriction standard persistent excitation suffices to ensure
exponentially fast convergence of the estimates to the actual
values of unknown parameters. Subsequently, our approach is
extended to cases in which the monotonicity assumption holds only
locally. We show that excitation with high-frequency of
oscillations is sufficient to ensure convergence. Two practically
relevant examples are given in order to illustrate the
effectiveness of the approach.
\end{abstract}


{ {\bf Keywords:} nonlinear parametrization, parameter estimation,
persistent excitation, exponential convergence, monotonic
functions}

 \vspace{10mm}
{ \small
{\bf Corresponding author:}{  \\Ivan Tyukin \\
                              Laboratory for Perceptual Dynamics,\\
                              RIKEN Brain Science Institute,\\
                              2-1, Hirosawa, Wako-shi, Saitama,\\
                              351-0198, Japan\\
                              phone: +81-48-462-1111 extension 7436\\
                              fax:   +81-48-467-7236\\
                              e-mail: tyukinivan@brain.riken.jp}
}
\section{Introduction}

Broad areas of applied and fundamental science require parameter
identification of nonlinear systems. Research of these systems has
made substantial progress in identification of both static and
dynamic linearly parameterized systems
\cite{Ljung_99},\cite{Eykhoff}, \cite{Bastin92},  as well as
static nonlinear ones
\cite{Box},\cite{Hansen92},\cite{Wilde},\cite{Kirkpatrick83}.
Parameter estimation of dynamic systems with nonlinear
parametrization, however, has long remained an open issue. One
possible way, in principle, to address this problem is to try to
derive the estimator for the general nonlinear case. In this hard
problem a breakthrough has resulted in an advanced method
\cite{Cao_2003}. This method applies to a large class of nonlinear
systems. Despite this major advancement the price for such
generality are several theoretical and practical limitations.
First, it is required that uncertainty be Lipshitz in time.
Second, extra control is needed in order to dominate the
nonlinearity during identification. The third and most important
restriction is the necessity to satisfy the {\it nonlinear
persistent excitation condition}, which computationally is more
difficult to check than conventional persistent excitation
assumptions \cite{Morgan_77}, \cite{Narendra89}\footnote{See also
\cite{Panteley_2001}, \cite{Zhang_96}, where relaxed formulations
of persistent excitation conditions are discussed.} and often is
not easy to satisfy if state-dependent nonlinearities are allowed.

An alternative strategy would be to consider a class of nonlinear
parameterizations that is narrower, but still sufficiently broad
to be practically relevant. Currently there are a number of
models, for instance Hammerstain (Wiener) models
\cite{Narendra_66},\cite{Pawlak_91},\cite{Garulli_2002},\cite{Bai_2003},
with specific restrictions on the nonlinearity in the parameters
that allow to avoid  the problems arising in the general
parametrization case. These models, however,  handle only static
input (output) nonlinearities. Local linear (nonlinear) model
techniques \cite{Johansen_95}, \cite{Verdult_2002},
\cite{Enqvist_CDC2002} constitute another promising tool. These
models, on the other hand, are not always physically plausible. In
practice in order to identify parameters of actual physical
processes in a system, it is often necessary with these models to
refit the data to the original nonlinearly parameterized model.

In this article we address the problem of parameter estimation of
dynamic nonlinear parameterized systems in a way that compromises
between the pros and cons of both strategies mentioned above. In
particular, we allow nonlinear state-dependent parametrization in
the model, while restricting the nonlinearities in parameters to
be of a certain practically relevant class. As a result we obtain
parameter estimation procedures that are not limited to Lipshitz
nonlinearities in time. These procedures neither require
domination of the nonlinearity, nor do they rely on nonlinear
persistent excitation conditions. On the other hand, the class of
nonlinear parameterizations that we propose is wide enough to
include a variety of models in physics, mechanics, physiology and
neural computation \cite{Armstrong_1993}, \cite{Pacejka91},
\cite{Boskovic_1995}, \cite{Abbott_2001}. For this new class of
parameterizations we show that conventional persistent excitation
conditions guarantee exponential convergence of the estimates to
the actual values of the parameters. Moreover, in case our
assumptions are satisfied only locally, sufficiently high
frequency of excitation still ensures convergence.



The paper is organized as follows. In Section 2 we formulate the
problem, Section 3 contains the main results of the paper, in
Section 3 we provide two practically relevant illustrative
applications of our method, and Section 4 concludes the paper.


\section{Problem Formulation}

Let the following system be given:
\begin{eqnarray}\label{system1}
\dot{\bfx}_1&=&\bff_1(\bfx)+\bfg_1(\bfx)u,\nonumber \\
\dot{\bfx}_2&=&\bff_2(\bfx,\thetavec)+\bfg_2(\bfx)u,
\end{eqnarray}
where
\[
\bfx_1=(x_{11},\dots,x_{1 m_1})^T\in \Real^{m_1}
\]
\[
\bfx_2=(x_{21},\dots,x_{2 m_2})^T\in \Real^{m_2}
\]
\[
\bfx=(x_{11},\dots,x_{1m_1},x_{21},\dots,x_{2m_2})^T\in \Real^{n}
\]
$\thetavec\in \Omega_\theta\in \Real^d$ is a vector of unknown
parameters, $u$ is the control input, and functions
$\bff_1:\Real^{n}\rightarrow \Real^{m_1}$,
$\bff_2:\Real^{n}\times\Real^d\rightarrow
\Real^{m_2}$,$\bfg_1:\Real^{n}\rightarrow \Real^{m_1}$,
$\bfg_2:\Real^{n}\rightarrow\Real^{m_2}$ are locally
bounded\footnote{Function $\bff(\bfx): \Real^{n}\rightarrow
\Real^m$ is said to be locally bounded if for any
$\|\bfx\|<\delta$ there exists  constant $D(\delta)>0$ such that
the following holds: $\|\bff(\bfx)\|\leq D(\delta)$.}. Vector
$\bfx\in\Real^n$ is a state vector, and vectors $\bfx_1$, $\bfx_2$
are referred to as {\it uncertainty independent} and {\it
uncertainty dependent} partitions of $\bfx$ respectively. We
assume that $\Omega_\theta$ is bounded, and therefore without loss
of generality it is safe to assume that $\Omega_\theta$ is a
closed ball or hypercube in $\Real^d$.

For the sake of compactness we introduce the following alternative
description for (\ref{system1}):
\begin{eqnarray}\label{system}
\dot{\bfx}=\bff(\bfx,\thetavec)+\bfg(\bfx)u,
\end{eqnarray}
where
\[
\bfg(\bfx)=(g_{11}(\bfx),\dots,g_{1m_1}(\bfx),g_{21}(\bfx),\dots,g_{2
m_2}(\bfx))^{T}
\]
\[
\bff(\bfx)=(f_{11}(\bfx),\dots,f_{1m_1}(\bfx),f_{21}(\bfx,\thetavec),\dots,f_{2
m_2}(\bfx,\thetavec))^{T}
\]

Our goal is to derive both the control function $u(\bfx,t)$ and
estimator $\hat{\thetavec}(t)$ such that all trajectories of the
system are bounded and the estimate $\hat{\thetavec}(t)$ converges
to unknown $\thetavec\in \Omega_\theta$ asymptotically. In
addition, in order to ensure boundedness of the trajectories we
will restrict all possible motions of system (\ref{system}) to an
admissible domain in the system state space. As a measure of
closeness of the trajectories to the desired solution we introduce
the smooth error function $\psi:\Real^n\times \Real\rightarrow
\Real, \ \psi\in C^1$. The function $\psi(\bfx,t)$ is bounded in
$t$ for every bounded $\bfx$. The target manifold, therefore,  is
given by
\[
\psi(\bfx,t)=0
\]
Consider the transverse dynamics of system (\ref{system}) with
respect to $\psi(\bfx,t)$:
\begin{eqnarray}\label{dpsi}
\dot{\psi}=L_{\bff(\bfx,\thetavec)}\psi(\bfx,t)+L_{\bfg(\bfx)}\psi(\bfx,t)u+\frac{\pd
\psi(\bfx,t)}{\pd t},
\end{eqnarray}
where $L_{\bff(\bfx,\thetavec)}$ is Lie derivative of function
$\psi(\bfx,t)$ with respect to vector field
$\bff(\bfx,\thetavec)$. Let us further assume that
$L_{\bfg(\bfx)}\psi(\bfx,t)$ is separated from zero (i.e. there
exists positive $\delta>0$ such that
$|L_{\bfg(\bfx)}\psi(\bfx,t)|>\delta$ for any $\bfx\in\Real^n$,
$t\in \Real_+$). This assumption automatically implies existence
of the inverse $L_{\bfg(\bfx)}\psi(\bfx,t)^{-1}$. Hence, we can
select control input $u$ from the following class of functions:
\begin{eqnarray}\label{control}
u(\bfx,\hat\thetavec,t)=(L_{\bfg(\bfx)}\psi(\bfx,t))^{-1}(-L_{\bff(\bfx,\hat{\thetavec})}\psi(\bfx,t)-\varphi(\psi)-\frac{\pd\psi(\bfx,t)}{\pd
t}),
\end{eqnarray}
where
\begin{eqnarray}\label{varphi}
\varphi:\Real\rightarrow \Real, \ \varphi(\psi)\in C^{1}, \
\varphi(\psi)\psi>0 \ \forall  \ \psi\neq 0, \ \
\lim_{\psi\rightarrow\infty}\int_0^{\psi}\varphi(\xi)d\xi=\infty.
\end{eqnarray}
Denoting
$L_{\bff(\bfx,\thetavec)}\psi(\bfx,t)=f(\bfx,\thetavec,t)$ and
taking into account (\ref{control}) we can rewrite equation
(\ref{dpsi}) in the following manner:
\begin{eqnarray}\label{error_model}
\dpsi=f(\bfx,\thetavec,t)-f(\bfx,\hat{\thetavec},t)-\varphi(\psi)
\end{eqnarray}

It is natural to require that boundedness of $\psi(\bfx,t)$
implies boundedness of the state insofar as $\psi(\bfx,t)$ stands
for the deviation from target manifold $\psi(\bfx,t)=0$. Let us
formally introduce this requirement in the following assumption:
\begin{assume}\label{assume:psi}  For the given function $\psi(\bfx,t)$ the following holds:
\[
\psi(\bfx,t)\in L_\infty \Rightarrow \bfx\in L_\infty
\]
\end{assume}
Assumption \ref{assume:psi} can be considered a bounded input -
bounded state assumption  for system (\ref{system1}) along the
constraint $\psi(\bfx,t)=\upsilon(t)$, where functions $u$ is
chosen to satisfy this requirement and signal $\upsilon(t)$ serves
as input. If, however, boundedness of the state is not required or
is achieved by extra control, Assumption \ref{assume:psi} can be
removed from the statements of our results or replaced, when
necessary, with the requirement for the function
$f(\bfx,\thetavec,t)$ to be globally bounded in $\bfx$ and locally
bounded in $\thetavec$.

So far the only deviation from standard descriptions of the
problem resides in the specification of function
$f(\bfx,\thetavec,t)$. Since a general parametrization of function
$f(\bfx,\thetavec,t)$ is methodologically difficult to deal with
but solutions provided for the restricted classes of
nonlinearities often yield physically implausible models, we have
opted to search for a new class of practically reasonable
parameterizations. Such a class should be able to include a
sufficiently broad range of physical models, in particular those
with nonlinear parametrization; they should also, in principle, be
able to handle arbitrary (in the class of smooth functions)
nonlinearity in states. As a candidate for such a parametrization
we suggest nonlinear functions that satisfy the following
assumption:
\begin{assume}[Monotonicity and Linear Growth Rate in Parameters]\label{assume:alpha}
There exists function $\alphavec(\bfx,t): \Real^{n}\times
\Real\rightarrow \Real^d$ and $D>0$ such that
\[
(f(\bfx,\hat{\thetavec},t)-f(\bfx,\thetavec,t))(\alphavec(\bfx,t)^{T}(\hat{\thetavec}-\thetavec))>0
\  \forall \ f(\bfx,\thetavec,t)\neq f(\bfx,\hat{\thetavec},t)
\]
\[|f(\bfx,\hat{\thetavec},t)-f(\bfx,\thetavec,t)|\leq
D |\alphavec(\bfx,t)^{T}(\hat{\thetavec}-\thetavec)|, \ D>0\]
\end{assume}
The first statement in Assumption \ref{assume:alpha} holds, for
example, for every smooth nonlinear function which is monotonic
with respect to a linear functional over a vector of parameters
$\thetavec$: $f(\bfx,\phivec(\bfx)^T\thetavec,t)$. The second
inequality is satisfied if the function
$f(\bfx,\phivec(\bfx)^T\thetavec,t)$ does not grow faster than a
linear function in variable $\phivec(\bfx)^T\thetavec$ for every
$\bfx\in \Real^n$. This set of conditions naturally extends
systems that are linear in parameters  to those with nonlinear
parametrization. In addition to linearly parameterized systems,
Assumption \ref{assume:alpha} covers a considerable large variety
of practically relevant models with  nonlinear parametrization.
These include effects of stiction forces \cite{Armstrong_1993},
slip and surface dependent friction given by the ``magic formula"
\cite{Pacejka91}, smooth saturation, and dead-zones in mechanical
systems.  The set of functions covered by Assumption
\ref{assume:alpha} further includes nonlinearities in models of
bio-reactors \cite{Boskovic_1995}. The class of functions
$f(\bfx,\thetavec,t)$ specified in Assumption \ref{assume:alpha}
can also serve as nonlinear replacement of the functions that are
linear in their parameters in a variety of piecewise approximation
models. Last but not least it includes sigmoid and Gaussian
nonlinearities, which are favored in neuro and fuzzy control and
mathematical models of neural processes \cite{Abbott_2001}.

In this article we attempt to resolve the following main issue:
how to design the estimator $\hat{\thetavec}(\bfx,t)$ which
ensures convergence of the estimates to the actual values of
a-priori unknown parameter $\thetavec$, and what further
restrictions (if any) on functions $f(\bfx,\thetavec,t)$ are to be
satisfied in order to guarantee such convergence?

\section{Main Results}



Let us introduce the following adaptation algorithm
\footnote{Parameter adjustment algorithms
(\ref{fin_forms_ours_tr1}) can be considered as generalizations of
the algorithms introduced earlier by the authors in
\cite{ECC_2003},\cite{t_fin_formsA&T2003}. In these works we
analyzed stabilizing properties of these algorithms in connection
with realizability issues. Parameter convergence and identifying
properties of algorithms (\ref{fin_forms_ours_tr1})
 were not addressed there.}:
\begin{eqnarray}\label{fin_forms_ours_tr1}
\hat{\thetavec}(\bfx,t)&=&\Gamma(\hat{\thetavec}_P(\bfx,t)+\hat{\thetavec}_I(t));\nonumber
\\ \hat{\thetavec}_P(\bfx,t)&=&
\psi(\bfx,t)\alphavec(\bfx,t)-\Psi(\bfx,t)\nonumber \\
\dot{\hat{\thetavec}}_I&=&\varphi(\psi(\bfx,t))\alphavec(\bfx,t)+{\pd
\Psi(\bfx,t)}/{\pd t}-\psi(\bfx,t)({\pd
\alphavec(\bfx,t)}/{\pd t})-\nonumber\\
& &  (\psi(\bfx,t)L_{\bff_1} \alphavec(\bfx,t)-L_{\bff_1}
\Psi(\bfx,t))-(\psi(\bfx,t)L_{\bfg_1}\alphavec(\bfx,t)-L_{\bfg_1}
\Psi(\bfx,t))u(\bfx,\hat{\thetavec},t)\nonumber \\
& & +
\beta(\bfx,t)(\bff_2(\bfx,\hat{\thetavec})+\bfg_2(\bfx)u(\bfx,\hat{\thetavec},t)),
\end{eqnarray}
where functions $\Psi(\bfx,t)$, $\beta(\bfx,t)$  satisfy the
following condition with respect to the vector-fields of system
(\ref{system1}) and function $\alpha(\bfx,t)$:
\begin{assume}\label{assume:explicit_realizability} There exists function $\Psi(\bfx,t)$ such that
\[
\frac{\pd \Psi(\bfx,t)}{\pd \bfx_2}-\psi(\bfx,t)\frac{\pd
\alpha(\bfx,t)}{\pd \bfx_2}=\beta(\bfx,t)
\]
where $\beta(\bfx,t)$ is ether zero or, if
$\bff_2(\bfx,\thetavec)$ is differentiable in $\thetavec$,
satisfies the following:
\[
\beta(\bfx,t)\mathcal{F(\bfx,\thetavec,\thetavec')}\geq 0 \ \
\forall \thetavec,\thetavec'\in \Omega_\theta, \ \bfx\in \Real^n
\]
\[
 \mathcal{F(\bfx,\thetavec,\thetavec')}=\int_0^1 \frac{\pd
\bff_2(\bfx, \bfs(\lambda))}{\pd \bfs} d\lambda, \ \
\bfs(\lambda)=\thetavec'\lambda+\thetavec(1-\lambda)
\]
\end{assume}
Assumption \ref{assume:explicit_realizability} can be viewed as a
kind of a structural restriction. Indeed, one can easily see that
it automatically holds for the cases where $\frac{\pd
\alphavec(\bfx,t)}{\pd \bfx_2}=0$, i.e. when function
$\alpha(\bfx,t)$ does not depend explicitly on vector $\bfx_2$,
which stands for the uncertainty-dependent partition of system
(\ref{system1}). Assumption \ref{assume:explicit_realizability}
holds also for one-dimensional uncertainty-dependent partitions if
function $\psi(\bfx,t)\frac{\pd \alphavec(\bfx,t)}{\bfx_2}$ is
Riemann-integrable with respect to  $\bfx_2$ (vector $\bfx_2$ is
one-dimensional in this case, $\beta(\bfx,t)=0$). Although it may
seem to be difficult  to find functions $\Psi(\bfx,t)$ satisfying
requirements of Assumption \ref{assume:explicit_realizability}, in
general the difficulty of the problem can significantly be reduced
by {\it embedding} the system dynamics into one of a higher order,
for which Assumption \ref{assume:explicit_realizability} is
satisfied a-priori. Sufficient conditions ensuring existence of
such embedding for the parameterizations of general structure are
provided in \cite{ECC_2003}. For systems where parametric
uncertainty  can be reduced to vector fields with low-triangular
structure the embedding is given in \cite{ALCOSP_2004}. An
alternative and the easiest way to construct this embedding is to
design a system of which the output $\hat{\bfx}_2(t)$ tracks
vector $\bfx_2$ with the prescribed level of performance by use of
high-gain robust observers. The former is then used in the
adjustment algorithm as replacement for the latter. This makes it
possible to reduce the problem either to one of already considered
cases of independence of $\alphavec(\bfx,t)$ on
uncertainty-dependent partitions or to single-dimension partitions
of $\bfx_2$. This technique is illustrated in detail in the
examples section.


Properties of system (\ref{system1}) with control (\ref{control})
and adaptation algorithm (\ref{fin_forms_ours_tr1}) are summarized
in Theorem \ref{stability_theorem} and Theorem
\ref{exp_stability_theorem}.

\begin{thm}[Stability and Convergence]\label{stability_theorem}
Let system (\ref{system1}), (\ref{control}),
(\ref{fin_forms_ours_tr1}) be given and Assumptions
\ref{assume:alpha}--\ref{assume:explicit_realizability} hold. Then

P1) $\varphi(\psi(t))\in L_2$, $\dpsi(t)\in L_2$;

P2) $\|\thetavec-\hat\thetavec(t)\|^{2}_{\Gamma^{-1}}$ is
non-increasing;

P3) $f((\bfx,\thetavec,t)-f(\bfx,\hat{\thetavec}(t),t))\in L_2$.

Furthermore,
\begin{eqnarray}\label{L_2_L_inf_performance}
 \|\varphi(\psi)\|_{2}^2 &\leq& 2
 Q(\psi)+\|\hat{\thetavec}(0)-\thetavec\|^{2}_{(2D\Gamma)^{-1}},
 \ \
 \|\dpsi\|_{2}^2 \leq  2 Q(\psi)+\|\hat{\thetavec}(0)-\thetavec\|^2_{(2D\Gamma)^{-1}} \nonumber \\
\|\psi\|_{\infty} &\leq& \Lambda
\left(Q(\psi)+\|\hat{\thetavec}(0)-\thetavec\|^2_{(4D\Gamma)^{-1}}\right),
\end{eqnarray}
where
$Q(\psi)=\int_{0}^{\psi(\bfx(0),0)}\varphi(\varsigma)d\varsigma$
and $\Lambda(d)=\max_{|\psi|}\{|\psi| \ | \
\int_{0}^{|\psi|}\varphi(\varsigma)d\varsigma=d\}$.

If Assumption \ref{assume:psi} is satisfied and function
$f(\bfx,\hat{\thetavec},t)$ is locally bounded with respect to
$\bfx$, $\hat{\thetavec}$ and uniformly bounded with respect to
$t$, then

P4) trajectories of the system are bounded and
$\psi(\bfx(t))\rightarrow 0$ as $t\rightarrow\infty$;

If in addition functions $\varphi, f(\bfx,\thetavec,t)\in C^1$,
derivative $ \pd {f(\bfx,\thetavec,t)}/{\pd t}$ is uniformly
bounded in $t$, function $\alphavec(\bfx,t)$ is locally bounded
with respect to $\bfx$ and uniformly bounded with respect to $t$,
then

P5) $\dpsi\rightarrow 0$ as  $t\rightarrow \infty$;
$f((\bfx,\thetavec,t)-f(\bfx,\hat{\thetavec}(t),t))\rightarrow 0$
 as  $t\rightarrow \infty$.

\end{thm}
Proofs of Theorem \ref{stability_theorem} and subsequent results
are given in the Appendix.

Theorem \ref{stability_theorem} ensures for algorithms
(\ref{fin_forms_ours_tr1}) asymptotic reaching of the control goal
and boundedness of the solutions of the closed-loop system. In
addition, it provides improved transient performance, which can be
characterized by a-priori computable $L_2$ norms for $\dpsi$ and
$\psi$. In the case where $f(\bfx,\thetavec,t)\neq
f(\bfx,\hat{\thetavec},t)$ along the system solutions it further
guarantees reduction of parametric uncertainties (property P2).


So far we have assumed that functions $\varphi(\psi)$  may vary
freely in the class of functions specified by condition
(\ref{varphi}). It is possible, however, to show that the
transient performance of system (\ref{system1}) with algorithms
(\ref{fin_forms_ours_tr1}) can further be improved when functions
$\varphi(\psi)$  are linear in $\psi$. An additional assumption on
the growth rate of function $f(\bfx,\thetavec)$  in $\thetavec$
will make the whole system exponentially stable. This new
assumption is formulated as follows:
\begin{assume}\label{assume:alpha_upper} For the given function
$f(\bfx,\thetavec)$ in (\ref{error_model}) and function
$\alpha(\bfx,t)$, satisfying Assumption \ref{assume:alpha}, there
exists a positive constant $D_1>0$  such that
\[|f(\bfx,\hat{\thetavec},t)-f(\bfx,\thetavec,t)|\geq
D_1 |\alphavec(\bfx,t)^{T}(\hat{\thetavec}-\thetavec)|
\]
\end{assume}
Assumption \ref{assume:alpha_upper} extends Assumption
\ref{assume:alpha} by stipulating a lower bound for the growth
rate of nonlinear function $f(\bfx,\thetavec)$ in $\thetavec$.
This assumption allows us to show that exponential convergence of
$\hat{\thetavec}$ to $\thetavec$ will automatically result in
exponentially fast convergence of function $\psi(\bfx,t)$ to the
origin. Furthermore, it ensures exponential convergence of
$\hat{\thetavec}$ to $\thetavec$ for any positive-definite
constant $\Gamma$. These results are formulated in the following
theorem:
\begin{thm}[Exponential Convergence]\label{exp_stability_theorem} Let Assumptions \ref{assume:alpha}--\ref{assume:explicit_realizability} hold and
$\varphi(\psi)=K \psi$, $K>0$. Then

P6) function $\psi(\bfx(t),t)$ converges exponentially fast into
the domain $|\psi(\bfx(t),t)|\leq
0.5\sqrt{\|\hat{\thetavec}(0)-\thetavec\|^2_{(KD\Gamma)^{-1}}}$.
Specifically, the following holds:   $|\psi(\bfx(t),t)|\leq
 |\psi(\bfx(0),0)| e^{-K t} + 0.5\sqrt{\|\hat{\thetavec}(0)-\thetavec\|^2_{(KD\Gamma)^{-1}}}$

Furthermore, let Assumption \ref{assume:psi} hold, function
$\alphavec(\bfx,t)$ be locally bounded with respect to $\bfx$ and
uniformly bounded in $t$; for any bounded $\bfx$ there exist
$D_1>0$ such
 that
$|f(\bfx,\hat{\thetavec},t)-f(\bfx,\thetavec,t)|\geq D_1
|\alphavec(\bfx,t)^{T}(\hat{\thetavec}-\thetavec)|$, function
$\alphavec(\bfx,t)$ is persistently exciting:
\begin{eqnarray}\label{PE}
\exists L>0, \ \delta>0: \
\int_{t}^{t+L}\alphavec(\bfx(\tau),\tau)\alphavec(\bfx(\tau),\tau)^{T}d\tau\geq
\delta I \ \forall t>0,
\end{eqnarray}
where $I\in R^{d\times d}$ -- identity matrix. Then

P7) both $\psi(\bfx(t),t)$ and $\|\hat{\thetavec}-\thetavec\|$
converge exponentially fast to the origin.
\end{thm}



It is desirable to notice that one can derive more precise
conditions for exponential convergence of the estimates with
algorithms (\ref{fin_forms_ours_tr1}) from the proofs of Theorems
\ref{stability_theorem}, \ref{exp_stability_theorem}. In
particular in the proofs we neglected term
$\beta(\bfx,t)\mathcal{F(\bfx,\thetavec,\hat{\thetavec})}$ in the
equations for derivatives $\dot{\hat{\thetavec}}$. The complete
set of conditions would be, therefore, as follows:
\begin{eqnarray}\label{PE_complete}
\exists L>0, \ \delta>0: & &
\int_{t}^{t+L}(\mathcal{F}_0^{T}(\bfx(t),\thetavec,\hat{\thetavec}(\tau),\tau)\mathcal{F}_0(\bfx(t),\thetavec,\hat{\thetavec}(\tau),\tau)
+\\\nonumber &
&\beta(\bfx,t)\mathcal{F}(\bfx,\thetavec,\hat{\thetavec}(\tau)))d\tau\geq
\delta I \ \forall t>0,
\end{eqnarray}
where matrix function
$\mathcal{F}(\bfx,\thetavec,\hat{\thetavec}(t))$ is defined as in
Assumption \ref{assume:explicit_realizability} and function $
\mathcal{F}_0(\bfx,\thetavec,\hat{\thetavec},t)$ is:
\[
 \mathcal{F}_0(\bfx,\thetavec,\hat{\thetavec},t)=\int_0^1 \frac{\pd
f(\bfx, \bfs(\lambda),t)}{\pd \bfs} d\lambda,  \ \
\bfs(\lambda)=\thetavec\lambda+\hat{\thetavec}(1-\lambda)
\]

So far we have shown that, for the class of nonlinearly
parameterized systems, there exist a control function and
parameter adjustment algorithms such that solutions of the whole
system are bounded, and parametric uncertainty is decreasing in
time.
We have shown also that in case of persistently excited functions
$\alphavec(\bfx,t)$ the estimates $\hat{\thetavec}(t)$ in
(\ref{fin_forms_ours_tr1}) converge exponentially fast to vector
$\thetavec$. These results will now be extended to a broader class
of nonlinearities. We replace Assumptions \ref{assume:alpha},
\ref{assume:alpha_upper} with their locally verified versions.
\begin{assume}\label{assume:interval_monotonic} For the given nonlinear function $f(\bfx,\thetavec)$ there exits
the following partition of the state space:
\[
\Omega_{\bfx}=\Omega_M(\bfx)\cup\Omega_{A}, \ \
\Omega_M(\bfx)=\bigcup_{j} \Omega_{M,j}(\bfx), \
\Omega_{A}=\Omega_{\bfx}/\Omega_{M}(\bfx)
\]
where
$\Omega_{M,j}(\bfx)=\{\bfx|(\bfx-\bfc_j)^{T}(\bfx-\bfc_j)\leq
r^2_j\}$ are the balls in $\Real^{n}$ where Assumptions
\ref{assume:alpha}, \ref{assume:alpha_upper} are satisfied for
every $\thetavec\in\Omega_\theta$ and corresponding functions
$\alphavec_j(\bfx,t)$ and constants $D_j$, $D_{1,j}$.
\end{assume}
A typical example of a nonlinear function which satisfies this
assumption is $\sin(\theta x)$, where the unknown parameter
$\theta$ belongs to a bounded interval. Another example is
$x^{\theta}$, $\theta\in[0,\infty)$. The last parametrization is
widely used in modelling physical ``power low" phenomena in nature
(see, for example \cite{Wu_1990}, where this function models
effects of nonlinear damping in muscles).

In the sequel we will denote control functions (\ref{control})
associated with parameter adjustment algorithms in $\Omega_{M,j}$
by symbol $u_{0,j}(\bfx,t)$. Once Assumptions \ref{assume:alpha},
\ref{assume:alpha_upper} hold only locally, we can guarantee
convergence of the estimates only if the state belongs to
$\Omega_M(\bfx)$. Therefore, extra control effort is needed. In
order to specify the desired feedback acting in the domain
$\Omega_A$ we introduce the following assumption on system
(\ref{system1}) dynamics:
\begin{assume}\label{assume:steer_exists} For any $\bfx_0\in\Omega_\bfx$ there exists a control function $u_j(\bfx,t)$
that steers the state $\bfx$ of system (\ref{system1}) into the
neighborhood of $\bfc_j: \ \|\bfx-\bfc_j\|^2\leq \delta_j^2, \
\delta<r^2_j$ in finite time.
\end{assume}
It should be noticed, however, that this assumption does not
require existence of stabilizing feedback, local or global in
Lyapunov sense, at the points $\bfx=\bfc_j$.

Let us finally consider the following control/identification
scheme:
\begin{eqnarray}\label{control_intervals}
& & \sigma_j=\left\{
        \begin{array}{ll}
        1-\sigma_j, & \bfx=B_{j,\sigma_j}\\
        \sigma_j, & \bfx\neq B_{j,\sigma_j}
        \end{array}
       \right. , \ \sigma_j(0)=\left\{
                                     \begin{array}{ll}
                                     0, & \|\bfx(0)-\bfc_j\|>\delta_j \\
                                     1, & \|\bfx(0)-\bfc_j\|\leq
                                     \delta_j
                                     \end{array}
                                     \right. \nonumber \\
& &  B_{j,0}=\{\bfx:\|\bfx-\bfc_j\|=\delta_j\}, \  B_{j,1}=\{\bfx:\|\bfx-\bfc_j\|=r_j\} \nonumber \\
u&=& (1-\sigma_j)u_j(\bfx,t) + \sigma_j u_{0,j}(\bfx,t) \nonumber \\
\hat{\thetavec}&=& \sigma_j\Gamma
(\hat{\thetavec}_P(\bfx,t)+\hat{\thetavec}_I(t)+C_j(t));\nonumber
\\ \hat{\thetavec}_P(\bfx,t)&=&\psi(\bfx,t)\alphavec_j(\bfx,t)-\Psi_j(\bfx,t)\nonumber \\
\dot{\hat{\thetavec}}_I&=&\sigma_j(\varphi(\psi(\bfx,t))\alphavec_j(\bfx,t)+{\pd
\Psi_j(\bfx,t)}/{\pd t}-\psi(\bfx,t)({\pd
\alphavec_j(\bfx,t)}/{\pd t})-\nonumber\\
& &  (\psi(\bfx,t)L_{\bff_1} \alphavec_j(\bfx,t)-L_{\bff_1}
\Psi_j(\bfx,t))-(\psi(\bfx,t)L_{\bfg_1}\alphavec_j(\bfx,t)-L_{\bfg_1}
\Psi_j(\bfx,t))u(\bfx,\hat{\thetavec},t)\nonumber \\
& & +
\beta_j(\bfx,t)(\bff_2(\bfx,\hat{\thetavec})+\bfg_2(\bfx)u(\bfx,\hat{\thetavec},t)))
\end{eqnarray}
\[
C_j(t)=(\thetavec_P(\bfx(t_{i-1}'),t_{i-1}')-\thetavec_P(\bfx(t_i),t_i)+C_j(t_{i-1}')),
\]
where $t_i$ are the time instants when $\bfx$ hits the domain
$\|\bfx-\bfc_j\|=\delta_j$ for $\sigma_j=1$ and $t_i'>t_i$ stands
for the time moments when the state $\bfx$ reaches
$\|\bfx-\bfc_j\|=r_j$ (for $\sigma_j=1$). Algorithm
(\ref{control_intervals}) includes algorithm
(\ref{fin_forms_ours_tr1}) as a part. It also includes switching
algorithm which specifies the time when parameter estimation
procedure (\ref{fin_forms_ours_tr1}) shall be turned ``on"/``off".
The identifying properties of this new algorithm follow from
Theorems \ref{stability_theorem}, \ref{exp_stability_theorem} and
are formulated in the following corollary:

\begin{cor}\label{cor:intervals} Let Assumptions \ref{assume:psi},
\ref{assume:interval_monotonic} hold and there exist at least one
$\alphavec_{j}(\bfx,t)$ such that Assumption
\ref{assume:explicit_realizability} is satisfied.  Then system
(\ref{system1}), (\ref{control_intervals}) trajectories are
bounded. If, in addition, function $\alphavec_j(\bfx,t)$ is
persistently exciting:
\begin{eqnarray}\label{PE1}
\exists \delta>0: \
\int_{t}^{t+L}\alphavec(\bfx(\tau),\tau)\alphavec(\bfx(\tau),\tau)^{T}d\tau\geq
\delta I \ \forall t>0,
\end{eqnarray}
with sufficiently small $L>0$, then parameters
$\hat{\thetavec}(t)$ converge to $\thetavec$ $t\rightarrow\infty$
monotonically with respect to the norm
$\|\hat{\thetavec}(t)-\thetavec\|$.
\end{cor}

A consequence of Corollary \ref{cor:intervals} is that increase of
excitation in functions $\alphavec_j(\bfx,t)$ results in an
extension of the class of nonlinearities suitable for our
approach. This is consistent with previously reported results
\cite{Cao_2003} on  parameter convergence in nonlinearly
parameterized systems. Whether extension of the class of
nonlinearities to more general functions renders it necessary to
increase excitation, however, is still an open issue\footnote{An
example is constructed in \cite{Cao_2003}, where nonlinear
persistent excitation condition holds for the given
parametrization, while the linear persistent excitation condition
for linear parametrization with respect to the same
parameter-independent function is not satisfied.}.

For illustration consider the following system as an application
of Corollary \ref{cor:intervals}:
\begin{eqnarray}\label{ex:sine}
\dot{x}_1&=&x_2\nonumber\\
\dot{x}_2&=&\sin(\theta x_1)+u,
\end{eqnarray}
where parameter $\theta\in \Omega_\theta=[0.6,1.4]$ is unknown
a-priori. The goal is to design  input $u(x_1,x_2,t)$ and
estimator $\hat{\theta}(t)$ such that  trajectories of the system
are bounded and $\hat{\theta}(t)\rightarrow\theta$ as
$t\rightarrow\infty$. For the given bounds of $\Omega_\theta$ we
first find the domain $\Omega_M$, satisfying Assumption
\ref{assume:interval_monotonic}:
\[
\Omega_M(\bfx)= \{\bfx \ | \ x_1\in [-3.38,-2.59]\}\cup \{\bfx \ |
\ x_1\in [-1.14,-1.14]\}\cup \{\bfx \ | \ x_1\in [2.59,3.38]\} =
\Omega_{M,1}\cup \Omega_{M,2}\cup \Omega_{M,3}
\]
Let us suppose that initial conditions of system (\ref{ex:sine})
are located most closely to the subset $\Omega_{M,1}(\bfx)=\{\bfx
\ | \ x_1\in [-3.38,-2.59]\}$. Then it is natural to assume that
the desired position of the plant for the purpose of
identification is in the center of $\Omega_{M,1}:
x_1=x_1^\ast=-2.985$. Function $\psi(\bfx,t)$ satisfying
Assumption \ref{assume:psi} is chosen in the following manner:
$\psi(x_1,x_2)=x_1+x_2-x_1^\ast$. Hence, according to
(\ref{control}), control inputs $u_{0,1}(\bfx,t)$, $u_{1}(\bfx,t)$
are given by equations:
\[
u_{0,1}(\bfx)=-x_2 - \sin(\hat{\theta},x_1)- \psi(x_1,x_2), \ \
u_{1}(\bfx)=-x_2 - \psi(x_1,x_2) - \sign(\psi(x_1,x_2))
\]
Function $\alpha(x_1,x_2)=-x_1$, and function $\hat{\theta}_P$,
 in (\ref{control_intervals}) is as follows:
\[
\hat{\theta}_P(\bfx)= \psi(x_1,x_2)\alpha(x_1,x_2)-\Psi(x_1,x_2),
\ \Psi(x_1,x_2)=(x_1-x_1^{\ast})x_2+ \frac{x_2^2}{2}.
\]
Function $\hat{\theta}_I$ now follows explicitly from equation
(\ref{control_intervals}).

In the next section we illustrate the application to and main
steps in the design of our algorithms for two diverse, challenging
and practically relevant problems. In the first example we apply
our approach to the optimal slip identification problem in brake
control systems. The second example provides a system for adaptive
content-dependent filtering and classification of visual
information.

\section{Examples}


{\it Example 1. Braking wheel control problem.} Consider the
problem of minimizing  the braking distance for a single wheel
rolling along a surface. The surface properties can vary depending
on the current position of the wheel. The wheel dynamics can be
given by the following system of differential equations
\cite{I_Petersen_2003}:
\begin{eqnarray}\label{ex:slip_model}
\dot{x}_1&=&-\frac{1}{m}F_s(F_n,\bfx,\theta)\nonumber\\
\dot{x}_2&=&\frac{1}{J}(F_s(F_n,\bfx,\theta)r-u)\nonumber \\
\dot{x}_3&=&-\frac{1}{x_1}((\frac{1}{m}(1-x_3)+\frac{r^2}{J})F_s(\bfx,\theta)-\frac{r}{J}u),
\end{eqnarray}
$x_1$ is longitudinal velocity, $x_2$ is angular velocity,
$x_3=(x_1-x_2)/x_1$ is wheel slip, $m$ is mass of the wheel, $J$
is moment of inertia, $r$ is radius of the wheel, $u$ is control
input (brake torque), $F_s(F_n,\bfx,\theta)$ is a function
specifying the tire-road friction force depending on the
surface-dependent parameter $\theta$ and the load force $F_n$.
This function, for example, can be derived from steady-state
behavior of the LuGre tire-road friction model
\cite{Canudas_1999}:
\[
F_s(F_n,\bfx,\theta)=F_n \sign(x_2)\frac{\frac{\sigma_0}{L}
g(x_2,x_3,\theta)\frac{x_3}{1-x_3}}{\frac{\sigma_0}{L}
\frac{x_3}{1-x_3}+g(x_2,x_3,\theta)},
\]
\[
g(x_2,x_3,\theta)=\theta(\mu_C + (\mu_{S}-\mu_{C})e^{-\frac{|r x_2
x_3|}{|1-x_3|v_s}}), \ \
\]
where $\mu_C$, $\mu_S$ are Coulomb and static  friction
coefficients, $v_s$ is the Stribeck velocity, $\sigma_0$ is the
normalized rubber longitudinal stiffness, $L$ is the length of the
road contact patch. In order to avoid singularities we assume, as
suggested in \cite{I_Petersen_2003}, that the system is turned off
when velocity $x_1$ reaches a small neighborhood of zero (in our
example we stopped simulations as soon as $x_1$ becomes less than
$5$ m/sec).

While the majority of the model parameters can be estimated
a-priori, the tire-road parameter $\theta$ is dependent on the
properties of the road surface. Therefore, on-line identification
of the parameter $\theta$ is desirable in order to compute the
optimal slip value
\begin{eqnarray}\label{ex:slip_value}
x_3^{\ast}=\arg\max_{x_3}{F_s(F_n,\bfx,\theta)}
\end{eqnarray}
which ensures the maximum deceleration force and therefore results
in the shortest braking distance.

The main loop controller is derived in accordance with the
standard certainty-equivalence principle and can be written as
follows:
\[
u(\bfx,\hat{\theta},x_3^\ast)=\frac{J}{r}((\frac{1}{m}(1-x_3)+\frac{r^2}{J})F_s(F_n,\bfx,\hat{\theta})-K_s
x_1(x_3-x_3^\ast)), \ K_s>0
\]
In order to estimate parameter $\theta$ by measuring the values of
variables $x_1,x_2$ and $x_3$, we construct the following
subsystem:
\[
\dot{\hat{x}}_3=-\frac{1}{x_1}((\frac{1}{m}(1-x_3)+\frac{r^2}{J})F_s(F_n,\bfx,\hat{\theta})-\frac{r}{J}u)+(x_3-\hat{x}_3)
\]
and consider dynamics of the error function
$\psi(\bfx,t)=\psi(x_3,\hat{x}_3)=x_3-\hat{x}_3$:
\begin{eqnarray}\label{ex:slip_error_model}
\dot{\psi}=-\psi+\frac{1}{x_1}((\frac{1}{m}(1-x_3)+\frac{r^2}{J})(F_s(F_n,\bfx,\theta)-F_s(F_n,\bfx,\hat{\theta}))
\end{eqnarray}
Function
$\frac{1}{x_1}((\frac{1}{m}(1-x_3)+\frac{r^2}{J})F_s(\bfx,\theta)$
is monotonic in $\theta$ and satisfies Assumptions
\ref{assume:alpha}, \ref{assume:alpha_upper} with
\[
\alpha(\bfx,t)=\frac{1}{x_1}((\frac{1}{m}(1-x_3)+\frac{r^2}{J})g(x_2,x_3,1)
\]
Therefore, in order to design an estimation scheme satisfying
assumptions of Theorem \ref{exp_stability_theorem} we shall find
functions $\Psi(\bfx,t)$, $\beta(\bfx,t)$ such that Assumption
\ref{assume:explicit_realizability} holds. Notice that every
equation in (\ref{ex:slip_model}) depends on unknown parameter
$\theta$ explicitly. Therefore, according to the introduced
terminology, there is no uncertainty independent partition of
system (\ref{ex:slip_model}), i.e. $\bfx=\bfx_2$. Let us choose
$\beta(\bfx,t)=0$. Then Assumption
\ref{assume:explicit_realizability} reduces to the following
equation:
\begin{eqnarray}\label{ex:explicit_realizability}
\frac{\pd \Psi(\bfx,t)}{\pd \bfx}=\psi(\bfx,t)\frac{\pd
\alpha(\bfx)}{\pd \bfx}
\end{eqnarray}
Instead of trying to solve this equation explicitly for function
$\Psi(\bfx,t)$  we {\it embed} system (\ref{ex:slip_model}),
(\ref{ex:slip_error_model}), as suggested in \cite{ECC_2003}, into
one of higher order, such that for the new set of equations
Assumption \ref{assume:explicit_realizability} will be reduced to
a case where $\alpha(\bfx,t)$ does not depend explicitly on
$\bfx_2$. In fact, the problem with Assumption
\ref{assume:explicit_realizability} would be solved if we replace
$\alpha(\bfx,t)$ with the function of time $\xi(t): \
\Real_+\rightarrow \Real$, of which the derivative is known. Let
us derive the required function $\xi(t)$.  Notice that function
$\alpha(\bfx,t)$ is continuous in its arguments and, moreover,
differentiable for $x_1\geq 5$. Therefore, given that the
right-hand side of  system (\ref{ex:slip_model}) is locally
bounded, we can conclude that function $\alpha(\bfx,t)$ has a
bounded derivative for bounded $\bfx$. The state, moreover, is
bounded as longitudinal velocity $x_1$ and angular velocity $x_2$
are bounded during the braking/acceleration regime, and relative
slip $x_3$ is bounded by the way it is defined in
(\ref{ex:slip_model}). Therefore, it is possible to track signal
$\alpha(\bfx(t),t)$ with arbitrary high precision by use of smooth
high-gain estimators. If estimators with discontinuous right-hand
sides are allowed then it is possible to provide exact tracking of
$\alpha(\bfx,t)$. Let us consider the following candidate for the
estimator of $\alpha(\bfx,t)$:
\begin{eqnarray}\label{ex:ext}
\dot{\xi}=-\frac{\pd \alpha}{\pd x_2}\frac{1}{J}u+\frac{\pd
\alpha}{\pd
x_3}\frac{1}{x_1}\frac{r}{J}u-K_\xi\varphi_\xi(\xi-\alpha(\bfx,t)),
\end{eqnarray}
where $K_\xi>0$ and
$\varphi_\xi(\xi-\alpha(\bfx,t))(\xi-\alpha(\bfx,t))\geq 0$ are to
be chosen to dominate the following sum
\[
\left(-\frac{\pd \alpha}{\pd x_1}\frac{1}{m}+\frac{\pd \alpha}{\pd
x_2}\frac{r}{J} - \frac{\pd \alpha}{\pd
x_3}\frac{1}{x_1}((\frac{1}{m}(1-x_3)+\frac{r^2}{J})\right)F_s(F_n,\bfx,\theta)
\]
for $x_1\in[5,40]$, $x_2\in[100,1]$, $x_3\in(0,1)$,
$\theta\in(0,2]$ and $|\xi-\alpha(\bfx,t)|>\varepsilon_0=0.001$.
Let us assume for the moment that function
$\alpha(\bfx,t)=\xi(t)$. Taking this property into account we can
extend system (\ref{ex:slip_error_model}) with equation
(\ref{ex:ext}) and replace function $\alpha(\bfx,t)$ in
description (\ref{fin_forms_ours_tr1}) of the algorithm with
function $\xi(t)$ of which the derivative is known. Hence, for the
new system Assumption \ref{assume:explicit_realizability} will be
automatically satisfied. Then according to
(\ref{fin_forms_ours_tr1}) and (\ref{ex:slip_error_model})
parameter adjustment algorithm will be given by the following
system:
\begin{eqnarray}\label{ex:slip_alg}
\hat{\theta}=-\gamma((x_3-\hat{x}_3)\xi+\hat{\theta}_I), \ \
\gamma=100
\end{eqnarray}
\[
\dot{\hat{\theta}}_I=(x_3-\hat{x}_3)(\xi-\dot{\xi})
\]
The only difference between algorithm (\ref{ex:slip_alg}) and
those which follow from explicit analytical solution of
(\ref{ex:explicit_realizability}) is in the residual term
$\varepsilon(t)=\xi-\alpha(\bfx,t)$, which can be made arbitrary
small. On the other hand, according to Theorem
\ref{exp_stability_theorem} the original system (without replacing
$\alpha(\bfx,t)$ with $\xi(t)$) is exponentially stable, which in
turn guarantees convergence of the estimates $\hat{\theta}$ to the
actual values of $\theta$ with algorithm (\ref{ex:slip_alg}) if
$\varepsilon(t)$ is sufficiently small\footnote{In this particular
example in addition to the exponential stability argument one can
easily derive from differential equations for $\hat{\theta}$:
$\dot{\hat{\theta}}=-\gamma
(\dpsi+\psi)(\alpha(\bfx,t)+\varepsilon(t))$ from the proofs of
Theorems \ref{stability_theorem}, \ref{exp_stability_theorem}
that. This equation implies exponential convergence of
$\hat{\theta}$ to $\theta$, provided that
$\alpha(\bfx,t)-\varepsilon(t)>\delta>0$ for some positive
constant $\delta$}.

We simulated system (\ref{ex:slip_model}) -- (\ref{ex:slip_alg})
with the following setup of parameters and initial conditions:
$\sigma_0=200$, $L=0.25$, $\mu_C=0.5$, $\mu_S=0.9$, $v_s=12.5$,
$r=0.3$, $m=200$, $J=0.23$, $F_n=3000$, $K_s=30$. The
effectiveness of estimation algorithm (\ref{ex:slip_alg}) could be
illustrated with Figure 1. Estimates $\hat{\theta}$ approach the
actual values of parameter $\theta$ sufficiently fast for the
controller to calculate the optimal slip value $x_3^\ast$ and
steer the system toward this point in real braking time.
Effectiveness of the proposed identification-based control can be
confirmed by comparing the braking distance in the system with
on-line estimation of $x_3^\ast$ according to formula
(\ref{ex:slip_value}) with $\theta=\hat{\theta}$ with the one, in
which the values of $x_3^\ast$ were kept constant (in the interval
$[0.1,0.2]$). For model parameters as presently given and road
condition given by the piece-wise constant function
\[
\theta(s)=\left\{
                    \begin{array}{ll}
                     0.3, & s\in[0,8]\\
                     1.3, & s\in(8,16]\\
                     0.7, & s\in(16,24]\\
                     0.4, & s\in(24,32]\\
                     1.5, & s\in(32,40]\\
                     0.6, & s\in(40,\infty]
                    \end{array}
                  \right., s=\int_0^{t}x_1(\tau)d\tau
\]
the simulated braking distance obtained with our on-line
estimation procedure of $x_3^\ast$ is $54.95$ meters. This result
compares favorably with the values obtained for preset values of
$x_3^\ast$, which range between $57.52$ and $55.32$ (for
$x_3^\ast=0.1$ and $x_3^\ast=0.2$ respectively).

\begin{figure}\label{fig:slip}
\begin{center}
\includegraphics[width=200pt]{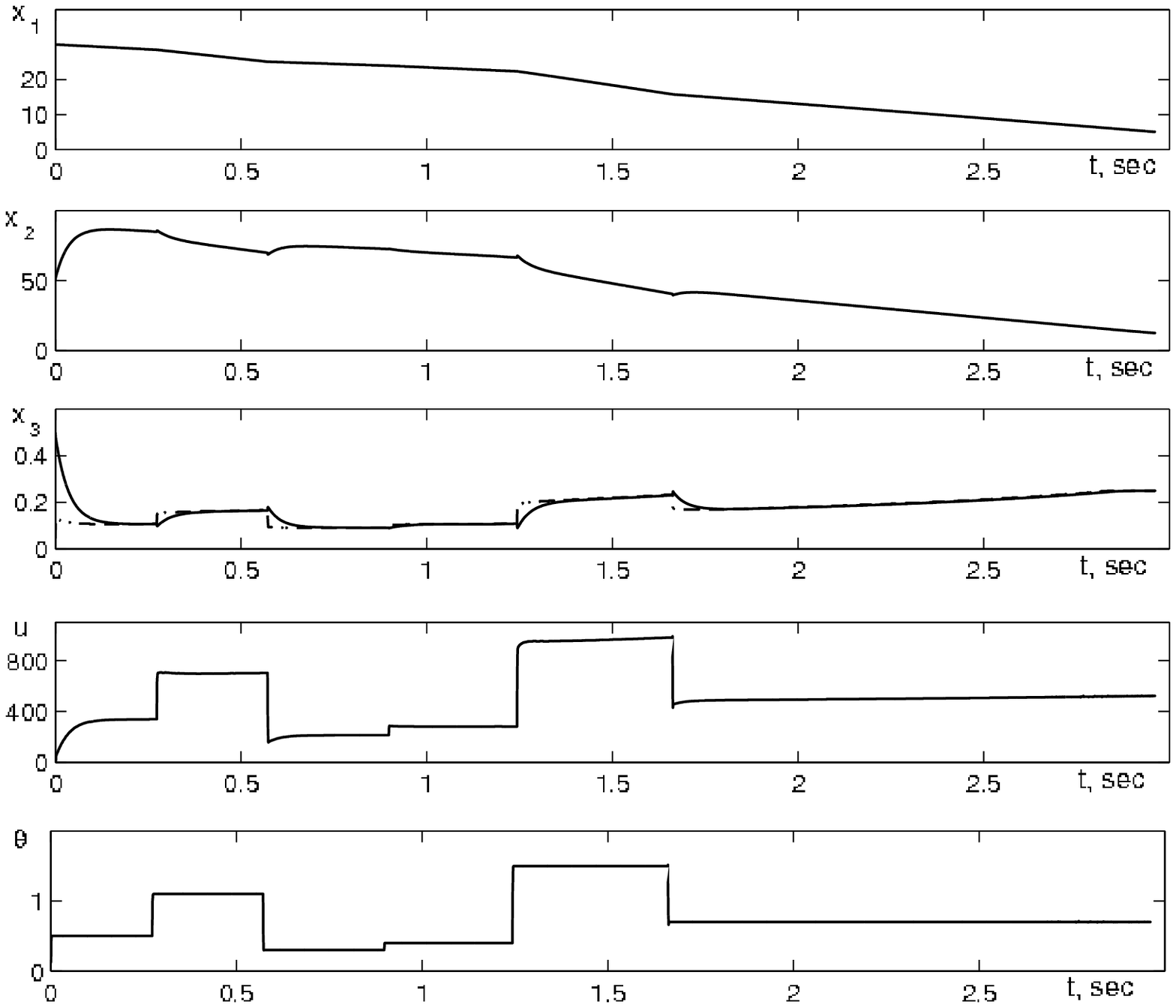}
\end{center}
\caption{Plots of the trajectories of system
(\ref{ex:slip_model})}
\end{figure}

{\it Example 2. Classification of occluded and linearly
nonseparable patterns}. Another illustrative application of our
parameter identification scheme is taken from the field of neural
computation and control of biological systems. In these domains
functions that are nonlinear in their parameters are widely used.
We will discuss an example involving a model of visual object
recognition system with adaptive identification of local spatial
features of the presented objects.

The problem is to identify two visual patterns given that they may
be out of focus (blurred), partially occluded by each other,
slightly distorted, or incomplete. A typical classification system
for this purpose consists of a two-layer neural network in which
the first, sensory layer feeds its outputs to a layer of decision
units.

We modeled a simple pattern recognition system, consisting of a
pattern input system and two image template systems. In the input
system, visual information arrives initially in a two-dimensional
array $(i, j)$ of sensors. The output of each sensor is mapped
onto the first layer, which consists of a two-dimensional array
$(i, j)$ of filters. The connectivity of the sensors to the
filters is one-to-all, but assures topographical projection by
means of connection weights. Maximum weight is given to
topographically corresponding units, neighboring ones receive
exponentially declining weights. The process is functionally
equivalent to spatial integration with an exponential kernel. The
output of these filters is projected topographically (one-to-one)
onto the second layer, which consists of a two-dimensional array
of decision units. Connections within each layer have not been
modeled at any of these levels. The architecture of the template
systems is identical to that of the input system, except that each
of the sensors is replaced by a binary value ("off" or "on"),
corresponding to the template of an image stored in memory. These
images can be represented by binary matrices $P_1$ and $P_2$
correspondingly.

The input and template systems are connected at the level of the
first and second layer. Connections at the level of the filters
layer are one-to-one, reciprocal but not symmetrical, between the
filters and their counterparts within each of the template
systems. Connections at the level of the decision layer are
all-to-all, reciprocal and symmetrical (for simplicity) between
the decision units and their counterparts in the template systems.
Also the two template systems are connected to each other in this
manner at this level.

Whereas the sensory units locally filter the spatial information
in the input, the decision nodes match it with the templates. The
information in the input system matches either of the templates,
to the degree that its units are synchronized with those of either
template. In addition, the location of the synchronized nodes in
the decision units indicates where the matching occurs in the
input.

The decision units and their counterparts in the template systems
can be modelled by the following ensemble of Hindmarsh and Rose
spiking neurons \cite{Hindmarsh_and_Rose}\footnote{For the sake of
compactness in the text hereafter  we omit indices $i,j$ in the
subscripts of the system variables}:
\begin{eqnarray}\label{ex: HR_model}
    \dot{x}_{1,k} &=& -a{x}_{1,k}^3+b{x}_{1,k}^2+x_{4,k}+x_{2,k}-x_{3,k}+u_k+I_0 \nonumber \\
    \dot{x}_{2,k} &=& c-d{x}_{1,k}^2-x_{2,k}\nonumber \\
    \dot{x}_{3,k} &=& \eps(s(x_{1,k}+x_0)-x_{3,k}), \nonumber \\
& & u_k=\gamma(x_{1,m}+x_{1,r}-2 x_{1,k}), \ \ m\neq r, k\neq m, \
\ k,m,r \in \{1,2,3\}
\end{eqnarray}
Parameters $a,b,c,d,s,x_0,\epsilon, I_0$ are all positive
constants with the following values: $a = 1$, $b = 3$, $c = 1$, $d
= 5$, $s = 4$, $x_0 = 1.6$, $\eps = 0.001$, $I_0=1.4$ as specified
by \cite{Hansel_1992}. Function $u_k=-2
x_{1,k}+\hat{x}_{1,k}+\bar{x}_{1,k}$ is the coupling function,
variables $\hat{x}_{1,k}$, $\bar{x}_{1,k}$ are the outputs of the
template systems at the level of the decision units and their
corresponding coupling functions are $\hat{u}_k=-2
\hat{x}_{1,k}+{x}_{1,k}+\bar{x}_{1,k}$, $\bar{u}_k=-2
\bar{x}_{1,k}+\hat{x}_{1,k}+{x}_{1,k}$.

Variable $x_{4,k}$ stands for input dependent current produced by
the sensory cell:
\begin{eqnarray}\label{ex:perception}
\dot{x}_{4,k}&=& \frac{1}{\tau}\left(\beta -x_{4,k} +
r(\theta_0,\bfs_k(t))
\right) \\
r(\theta_0,\bfs_k(t))&=&\sum_{i=1}^N\sum_{j=1}^{N}
e^{-\frac{|i(k)-i|+|j(k)-j|}{\theta_0}}
s_{k,i,j}\delta(t-\tau_{i,j})\nonumber, \ \beta>0
\end{eqnarray}
$i(k),j(k)$ specify the position of the $k$-th unit, $\tau>0$ is
the integration parameter, $s_{k,i,j}$ denotes intensity of the
$(i,j)$-th element in the image,
$\delta(t-\tau_{i,j}):\Real\rightarrow\Real_{\geq 0}$ stand for
the pattern-induced signals (impulses of unit amplitude and width
$\Delta=0.05 T$ at $t=\tau_{i,j}$, where $T$ is the period of
generation of each impulse). Exponential functions in
$r(\theta_0,\bfs_k(t))$  represent the distribution of the
weights. Numbers $\tau_{i,j}$ stand for time-delays in the
transmission of the signal from a sensor to the filters. This
delay is variable due to the difference in properties of the
transmission cables. It is given by the ration of cable length and
width, which is a simplification of actual signal transmission on
neural systems.  These delays together with the exponentially
decaying amplitudes of $s_{k,i,j}$ in space form the {\it
receptive field} of a filter. One of the main properties of such
an organization is that input signal $\bfs_k(t)$ is distributed in
time and space, providing in principle a unique spatiotemporal
signature for every different static visual pattern.

For the template systems the pattern-induced currents evolve
according to the following equations
\begin{eqnarray}\label{ex:adjusting_filters}
\dot{\hat{x}}_{4,k}&=& \frac{1}{\tau}\left( -\beta \hat{x}_{4,k} +
r({\theta_0},\hat{\bfs}_k(t,\theta_1)) \right),\nonumber \\
\dot{\bar{x}}_{4,k}&=& \frac{1}{\tau}\left( -\beta \bar{x}_{4,k} +
r({\theta_0},\bar{\bfs}_k(t,\theta_1)) \right)  \ \beta>0
\nonumber
\end{eqnarray}
where functions $\hat{\bfs}_k(t)$ and $\bar{\bfs}_k(t)$ are the
outputs of the spacial filters with parameter $\theta_1$:
\[
\hat{s}_{k,i,j}(t)=\delta(t-\tau_{i,j})\sum_{m=1}^{N}\sum_{r=1}^N
e^{\frac{-|i-m|-|r-j|}{\theta_1}}P_{1,m,r}, \ \
\bar{s}_{k,i,j}(t)=\delta(t-\tau_{i,j})\sum_{m=1}^{N}\sum_{r=1}^N
e^{\frac{-|i-m|-|r-j|}{\theta_1}}P_{2,m,r}
\]
This filters model effects of changes in intensity of the light,
focal adaptation and sharpness of the templates in presented
visual patterns.

The problem for such architecture is the following: if the picture
is not stable in time and perturbed by unmeasured changes in
focus, then how can the decision units reach detectable synchrony?

Technically, the solution would be to adjust parameters $\theta_1$
in (\ref{ex:adjusting_filters}) in response to distortion in the
input patters. The difficulty, however, is that the functions
$r(\theta,\hat{\bfs}_k(t,\theta_1))$
$r(\theta,\bar{\bfs}_k(t,\theta_1))$ are nonlinear in parameter
$\theta_1$. Classical linear identification schemes result in a
prohibitively large dimension of the estimator (in our simplified
example a $100\times 100$ sensory field will require $10^{10}$
independent parameters in each cell). A further problem is that of
performance in terms of robust identification of the parameter for
slightly distorted patterns, if no persistent excitation in
$\bfs_k(t)$ is assumed.  For these reasons we would require new
adjustment algorithms to estimate parameters $\hat{\theta}_1$,
$\bar{\theta}_1$ in the junctions of nodes (\ref{ex: HR_model}),
(\ref{ex:perception}) in order to ensure adaptive properties of
the classifier together with practical realizability and
reliability.

In order to derive the estimation algorithms for the parameters
$\theta_1$ in the templates we introduce the following error
functions: $\hat{\psi}_k=x_{4,k}-\hat{x}_{4,k}$,
$\bar{\psi}_k=x_{4,k}-\bar{x}_{4,k}$. Dynamics of function
$\hat{\psi}_k(t)$, for example, follows from
(\ref{ex:adjusting_filters}) and (\ref{ex:perception}):
\begin{eqnarray}\label{ex:model_synch_rfields}
\dot{\hat\psi}=-\frac{\beta}{\tau}\hat{\psi}+(r(\theta_0,\bfs_k(t))-r({\theta}_0,\hat{\bfs}_k(t,\hat{\theta}_1)))
\end{eqnarray}
If the image contains the distorted template locally then the
following equality holds:
${\bfs}_k(t)=\hat{\bfs}_k(t,\hat{\theta}_1)$, and equation
(\ref{ex:model_synch_rfields}) becomes as follows:
\begin{eqnarray}\label{ex:model_synch_rfields2}
\dot{\hat\psi}=-\frac{\beta}{\tau}\hat{\psi}+(r(\theta_0,\hat{\bfs}_k(t,{\theta}_1))-r({\theta}_0,\hat{\bfs}_k(t,\hat{\theta}_1)))
\end{eqnarray}
Function $\hat{\bfs}_k(t,\hat{\theta}_1)$ is monotonic in
$\theta_1$ with respect to $\theta_1$ (i.e. with constant
$\alpha(\bfx,t)$). Hence, Assumptions \ref{assume:alpha},
\ref{assume:alpha_upper} are at least locally satisfied. Moreover,
function $\alpha(\cdot)$ does not depend on $\bfx$, so Assumption
\ref{assume:explicit_realizability} is satisfied as well.
Therefore, applying  Theorem \ref{exp_stability_theorem} we can
derive parameter adjustment algorithm in the following form
\begin{eqnarray}\label{ex:alg_rfield}
\hat{\theta}_1&=& x_{4,k}-\hat{x}_{4,k}+\hat{\theta}_{I,1}\nonumber \\
\dot{\hat{\theta}}_{I,1}&=& \frac{\beta}{\tau}
\left(x_{4,k}-\hat{x}_{4,k}\right)
\end{eqnarray}
According to Theorem \ref{exp_stability_theorem}, algorithm
(\ref{ex:alg_rfield}) combined  with
(\ref{ex:model_synch_rfields}) result in the desired estimator of
$\theta_1$, which in addition guarantees exponential stability of
the whole system with respect to the small perturbations in the
presented patterns.

To illustrate the performance of our classifier, a square and a
cross  (Figures 2 A and B) were used as reference patterns in the
template system. They were distorted and combined, one partially
occluding the other as in Figure 2 C. Blurred versions of Figure 2
C, shown in 2 D and E, were presented to the system.  The task was
to recognize the patterns at their corresponding locations. In our
simulations we used the following values of the model parameters:
$\beta=0.02$, $\tau=0.01$, $\tau_{i,j}$ were set chosen in the
interval from $0$ to $100$, $\gamma=1$, $x_{1,k}(0)=-1.6$,
$x_{2,k}(0)=-11.83$, $x_{3,k}(0)=1.46$, $I_0=1.4$, $x_{4,k}(0)=0$,
$\hat{\theta}_I(0)=1$.

Figure 3 shows the responses of two of the decision nodes. Those
decision nodes that are topographical projections from regions
where the square appeared are synchronized with their counterparts
in the "square" template system.  Likewise, those which correspond
to regions where the cross appeared were synchronzied with their
counterparts in the "cross" template system. The synchrony occurs
because the for the counterparts the parameters $\hat{\theta}_1$ (
or $\bar{\theta}_1$) converge to their true values.  Evolution of
the estimates of $\theta_1$ is shown in Figure 4.

Decision nodes at regions where no pattern was presented and their
counterparts in the template systems fail to reach synchrony. The
**theta parameters of these units, however, remain in a bounded
domain.

\begin{figure}\label{fig:stimulus_line}
\begin{center}
\includegraphics[width=200pt]{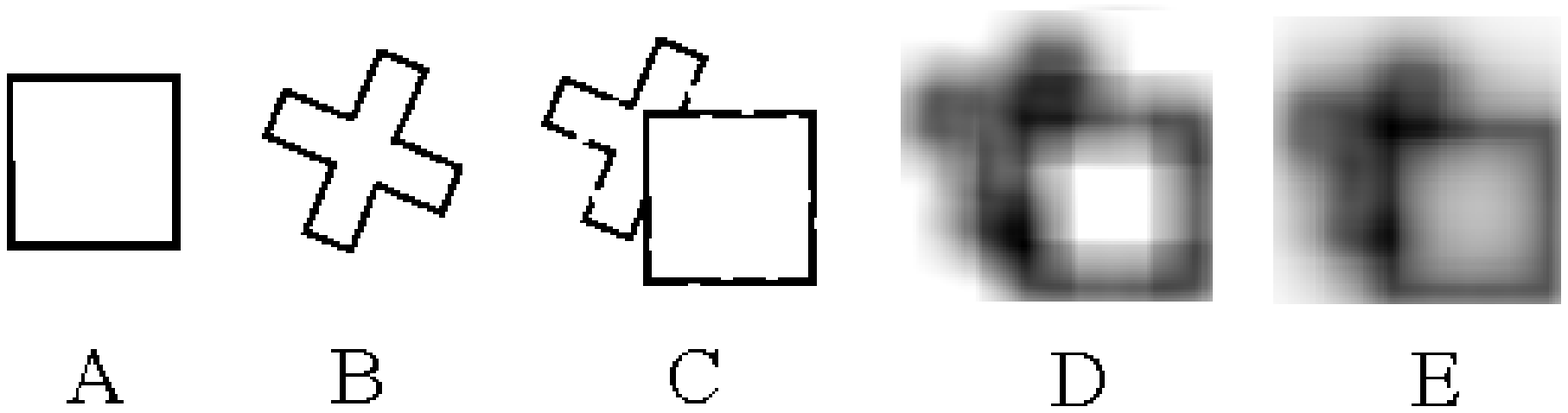}
\end{center}
\caption{Reference patterns of the template system (A,B);
distorted, combined pattern (C) ; blurred versions of the
distorted pattern (D,E)}
\end{figure}

\begin{figure}\label{fig:retina_response_distorted}
\begin{center}
\includegraphics[width=200pt]{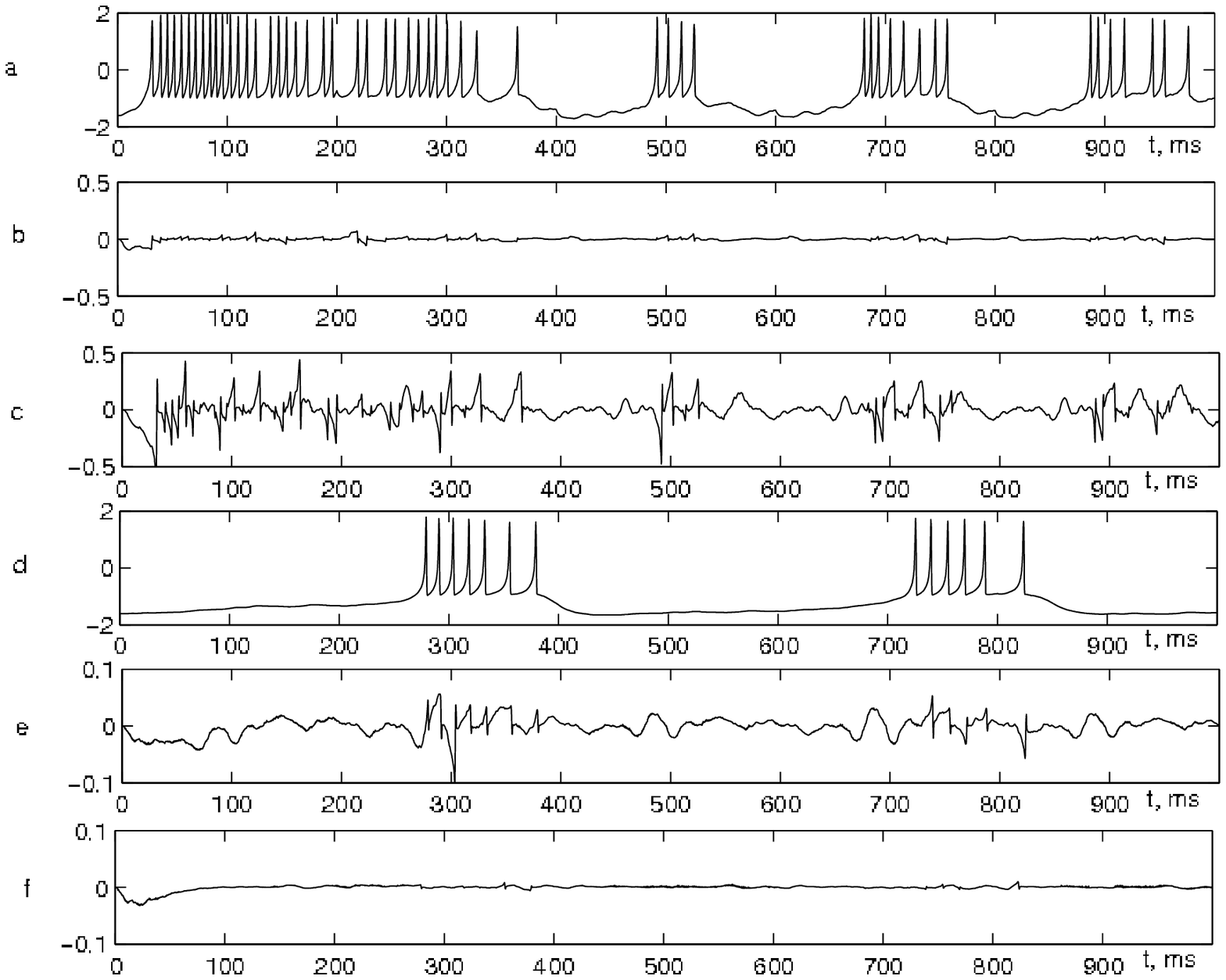}
\end{center}
\caption{Response of the decision nodes located in points
$(90,90)$ (plots $a$--$c$) and (5,20) (plots $d$--$e$)
respectively. Plots $a,d$ contain trajectories $x_{1,k}(t)$ of the
decision cells corresponding to the actual image, plots $b,e$
reflect differences in the responses between the actual scene and
memorized pattern "rectangle", plots $c,f$ show deviations in
perception of the scene and pattern "cross"}
\end{figure}

\begin{figure}\label{fig:retina_parameters}
\begin{center}
\includegraphics[width=200pt]{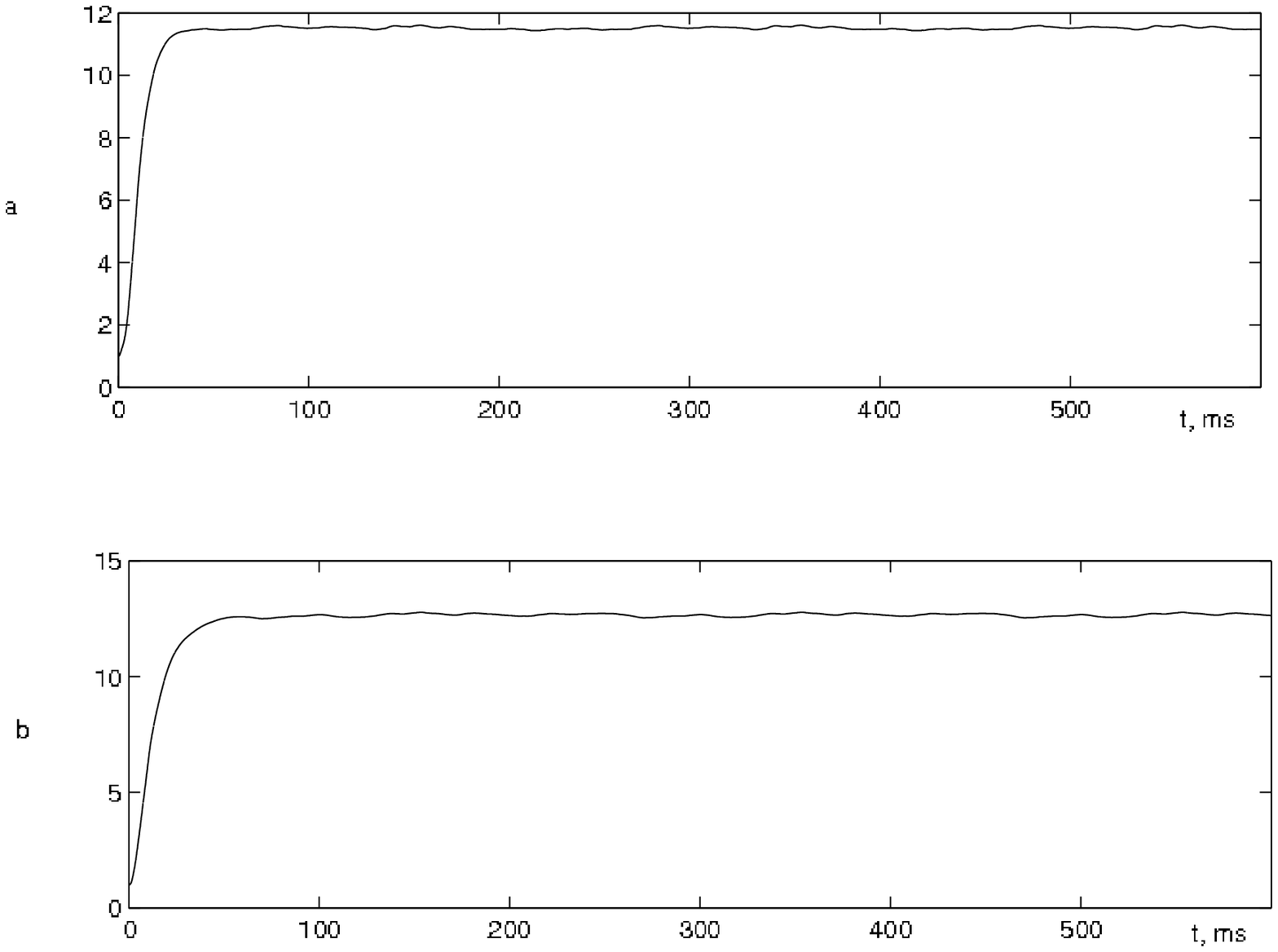}
\end{center}
\caption{Trajectories of $\hat{\theta}_1(t)$ in the decision cell
located at the point (90,90) (plot a), and $\bar{\theta}_1(t)$ in
the decision cell  located at the point $(5,20)$ (plot b).}
\end{figure}

\section{Conclusion}

We proposed a new class of parameterizations for nonlinearly
parameterized models. Instead of aiming at a general solution for
the problem of nonlinearity in the parameters, parametrization was
restricted to a set of smooth functions, which are monotonic with
respect to a linear functional in the parameters. For this new
class, estimation algorithms were introduced and analyzed. It was
been shown  that standard linear persistent excitation conditions
suffice to ensure exponentially fast convergence of the estimates
to the actual values of unknown parameters. If, however, the
monotonicity assumption holds only locally in the system state
space, excitation with sufficiently high-frequency of oscillations
is needed  to ensure convergence. It is also desirable to notice
that in case of linear parametrization the proposed parameter
estimation schemes  allow to estimate the unknowns in a dynamical
system without asking for usual filtered transformations, thus
reducing the number of integrators in the estimator.

Two rather distinct applications of our method were provided as
examples. One is devoted to on-line identification of the optimal
slip in a braking wheel. The second example touches on the problem
of dynamic recognition of visual patterns in artificial neural
networks. Both problems may be considered to have practical
significance. The effectiveness of the solution to these problems
leads us to expect that this method can successfully be
implemented in other applications.

\section{Appendix}

{\it Proof of Theorem \ref{stability_theorem}.} Let us first
calculate time-derivative  of function $\hat{\thetavec}(\bfx,t)$:
$\dot{\hat{\thetavec}}(\bfx,t)=\Gamma({\dot{\hat{\thetavec}}_{P}}+\dot{\hat\thetavec}_I)=\Gamma(\dpsi\alphavec(\bfx,t)+\psi\dot{\alphavec}(\bfx,t)-\dot{\Psi}(\bfx,t)+\dot{\hat\thetavec}_I)$.
Notice that
\begin{eqnarray}\label{t2_1}
&&\psi\dot{\alphavec}(\bfx,t)-\dot{\Psi}(\bfx,t)+\dot{\hat{\thetavec}}_I=\psi(\bfx,t)\frac{\pd
\alphavec(\bfx,t)}{\pd \bfx_1}\dot{\bfx}_1+\psi(\bfx,t)\frac{\pd
\alphavec(\bfx)}{\pd \bfx_2}\dot{\bfx}_2 + \psi(\bfx,t)\frac{\pd
\alphavec(\bfx,t)}{\pd t}-\nonumber\\
& & \frac{\pd \Psi(\bfx,t)}{\pd \bfx_1}\dot{\bfx}_1-\frac{\pd
\Psi(\bfx,t)}{\pd \bfx_2}\dot{\bfx}_2-\frac{\pd \Psi(\bfx,t)}{\pd
t}+\dot{\hat\thetavec}_I
\end{eqnarray}
According to Assumption \ref{assume:explicit_realizability},
$\frac{\pd \Psi(\bfx,t)}{\pd \bfx_2}=\psi(\bfx,t)\frac{\pd
\alphavec(\bfx,t)}{\pd \bfx_2}+\beta(\bfx,t)$. Then taking into
account (\ref{t2_1}), we can obtain
\begin{eqnarray}\label{t2_2}
& &
\psi\dot{\alphavec}(\bfx,t)-\dot{\Psi}(\bfx,t)+\dot{\hat{\thetavec}}_I=\left(\psi(\bfx,t)\frac{\pd
\alphavec(\bfx,t)}{\pd \bfx_1}-\frac{\pd \Psi}{\pd \bfx_1
}\right)\dot{\bfx}_1+\psi(\bfx,t)\frac{\pd \alphavec(\bfx,t)}{\pd
t}-\frac{\Psi(\bfx,t)}{\pd t}-\nonumber
\\
& &
\beta(\bfx,t)(\bff_2(\bfx,\thetavec)+\bfg_2(\bfx)u)+\dot{\hat{\thetavec}}_I
\end{eqnarray}
Notice that according to the proposed notation we can rewrite the
term $\left(\psi(\bfx,t)\frac{\pd \alphavec(\bfx,t)}{\pd
\bfx_1}-\frac{\pd \Psi}{\pd \bfx_1 }\right)\dot{\bfx}_1$ in the
following form: $\left(\psi(\bfx,t)L_{\bff_1}
\alphavec(\bfx,t)-L_{\bff_1} \Psi(\bfx,t)\right)+
\left(\psi(\bfx,t)L_{\bfg_1} \alphavec(\bfx,t)-L_{\bfg_1}
\Psi(\bfx,t)\right)u(\bfx,\hat{\thetavec},t)$. Hence it follows
from (\ref{fin_forms_ours_tr1}) and (\ref{t2_2}) that
$\psi\dot{\alphavec}(\bfx,t)-\dot{\Psi}(\bfx,t)+\dot{\hat{\thetavec}}_I=\varphi(\psi)\alphavec(\bfx,t)-\beta(\bfx,t)(\bff_2(\bfx,\thetavec)-\bff_2(\bfx,\hat{\thetavec}))$.
Therefore derivative $\dot{\hat\thetavec}(\bfx,t)$ can be written
in the following way:
\begin{eqnarray}\label{algorithm_dpsi}
\dot{\hat{\thetavec}}=\Gamma((\dpsi+\varphi(\psi))\alphavec(\bfx,t)-\beta(\bfx,t)(\bff_2(\bfx,\thetavec)-\bff_2(\bfx,\hat{\thetavec})))
\end{eqnarray}
Consider the following positive-definite function:
$V_{\hat{\thetavec}}(\hat{\thetavec},\thetavec)=
\frac{1}{2}\|\hat{\thetavec}-\thetavec\|^2_{\Gamma^{-1}}$. Its
time-derivative according to equations (\ref{algorithm_dpsi}) can
be obtained as follows:
\[
\dot{V}_{\hat{\thetavec}}(\hat{\thetavec},\thetavec)=(\varphi(\psi)+\dpsi)(\hat{\thetavec}-\thetavec)^{T}\alphavec(\bfx,t)-(\hat{\thetavec}-{\thetavec})^T\beta(\bfx,t)(\bff_2(\bfx,\thetavec)-\bff_2(\bfx,\hat{\thetavec}))
\]
Let $\beta(\bfx,t)\neq 0$, then  consider the following difference
$\bff_2(\bfx,\thetavec)-\bff_2(\bfx,\hat{\thetavec})$. Applying
Hadamard's lemma we represent this difference in the following
way:
\[
\bff_2(\bfx,\thetavec)-\bff_2(\bfx,\hat{\thetavec})= \int_0^1
\frac{\pd \bff_2(\bfx, \bfs(\lambda))}{\pd \bfs} d\lambda
(\thetavec-\hat{\thetavec}), \ \
\bfs(\lambda)=\thetavec\lambda+\hat{\thetavec}(1-\lambda)
\]
Therefore, according to Assumption
\ref{assume:explicit_realizability} function
$(\hat{\thetavec}-{\thetavec})^T\beta(\bfx,t)(\bff_2(\bfx,\thetavec)-\bff_2(\bfx,\hat{\thetavec}))$
is positive semi-definite, hence using Assumption
\ref{assume:alpha} and equality (\ref{dpsi}) we can estimate
derivative $\dot{V}_{\hat{\thetavec}}$ as follows
\begin{eqnarray}\label{parameric_deviation_derivative}
\dot{V}_{\hat{\thetavec}}(\hat{\thetavec},\thetavec)=-(f(\bfx,\hat{\thetavec},t)-f(\bfx,\thetavec,t))(\hat{\thetavec}-\thetavec)^{T}\alphavec(\bfx,t)
\leq-D(f(\bfx,\hat{\thetavec},t)-f(\bfx,\thetavec,t))^2=-D(\varphi(\psi)+\dpsi)^{2}\leq0
\end{eqnarray}
Therefore $V_{\hat{\thetavec}}$ is non-increasing (property P2) is
proven). Furthermore, integration of $\dot{V}_{\hat{\thetavec}}$
with respect to time results in
\[
V_{\hat{\thetavec}}(\hat{\thetavec}(0),\thetavec)-V_{\hat{\thetavec}}(\hat{\thetavec}(t),\thetavec)\geq
D\int_{0}^{t}(\dpsi(\tau)+\varphi(\psi(\tau)))^2 d\tau \geq 0.
\]
Function $V_{\hat{\thetavec}}$ is non-increasing and bounded from
below as  $V_{\hat{\thetavec}}\geq 0$, therefore
\[
D\int_{0}^{t}(\dpsi(\tau)+\varphi(\psi(\tau)))^2  d\tau \leq
V_{\hat{\thetavec}}(\hat{\thetavec}(0),\thetavec)<\infty.
\]
Hence
$(\varphi(\psi)+\dpsi)=(f(\bfx,\thetavec,t)-f(\bfx,\hat{\thetavec},t))=(f(\bfx,\thetavec,t)-f(\bfx,\hat{\thetavec},t))\in
L_2$ (property P3)).

To prove property P1) let us consider the following function:
$V(\psi,\hat{\thetavec},\thetavec)=2D Q(\psi)
+V_{\hat{\thetavec}}(\hat{\thetavec},\thetavec)$, where
$Q(\psi)=\int_0^{\psi}\varphi(\varsigma)d\varsigma$. Function
$V(\psi,\hat{\thetavec})$ is positive-definite with respect to
$\psi(\bfx,t)$ and $\hat{\thetavec}-\thetavec$. Its
time-derivative obeys inequality:
$\dot{V}(\psi,\hat{\thetavec},\thetavec)\leq2D\varphi(\psi)\dpsi -
D(\dpsi+\varphi(\psi))^{2}=-D\varphi^{2}(\psi)-D\dpsi^{2} \leq 0$.

Therefore, function $V(\psi,\hat{\thetavec},\thetavec)$ is bounded
and non-increasing. Furthermore
\begin{eqnarray}\label{t1_L2}
\infty>V(\psi(\bfx(0),0),\hat{\thetavec}(0),\thetavec)\geq
V(\psi(\bfx(0),0),\hat{\thetavec}(0),\thetavec)-V(\psi(\bfx(t),t),\hat{\thetavec}(t),\thetavec)&\geq&
D\int_{0}^t\varphi^2(\psi(\bfx(\tau),\tau))d\tau\geq 0\nonumber \\
\infty>V(\psi(\bfx(0),0),\hat{\thetavec}(0),\thetavec)\geq
V(\psi(\bfx(0),0),\hat{\thetavec}(0),\thetavec)-V(\psi(\bfx(t),t),\hat{\thetavec}(t),\thetavec)&\geq&
D\int_{0}^t\dpsi^2(\tau)d\tau\geq 0.
\end{eqnarray}
or, equivalently, $\dpsi(t)\in L_2$, $\varphi(\psi(t))\in L_2$.
Hence, property P1) is proven as well. The $L_2$ norm bounds
(\ref{L_2_L_inf_performance}) for $\varphi(\psi)$ and $\dpsi$
follow immediately from inequality (\ref{t1_L2}):
\begin{eqnarray}
\|\varphi(\psi)\|_2^2 & \leq&
{D}^{-1}V(\psi(\bfx(0),0),\hat{\thetavec}(0),\thetavec), \
\|\dpsi\|_2^2  \leq
{D}^{-1}V(\psi(\bfx(0),0),\hat{\thetavec}(0),\thetavec)\nonumber
\end{eqnarray}
The $L_\infty$ norm bound for $\psi(\bfx(t),t)$ results from the
inequality:
$V(\psi(\bfx(0),0),\hat{\thetavec}(0),\thetavec)-V(\psi(\bfx(t),t),\hat{\thetavec}(t),\thetavec)\geq0$.
Consider function $\Lambda$ defined as
$\Lambda(d)=\max_{|\psi|}\{|\psi| \ | \
\int_{0}^{|\psi|}\varphi(\varsigma)d\varsigma=d\}$ and notice that
it is monotonic and nondecreasing. Therefore, given that
$\int_{0}^{\psi(\bfx(t),t)}\varphi(\varsigma)d\varsigma\leq
\frac{1}{2D}V(\psi(\bfx(0),0),\hat{\thetavec}(0),\thetavec)$ we
can conclude that $|\psi|\leq
\Lambda\left(\frac{1}{2D}V(\psi(\bfx(0),0),\hat{\thetavec}(0),\thetavec)\right)$.
To prove property P4) notice that function
$V(\psi(\bfx(t),t),\hat{\thetavec}(t),\thetavec)$ is bounded.
Hence, as follows from condition (\ref{varphi}), function
$\psi(\bfx(t),t)$ is bounded as well. According to Assumption
\ref{assume:psi} boundedness of $\psi(\bfx(t),t)$ implies
boundedness of the state $\bfx$. In addition it is assumed that
$f(\bfx,\hat{\thetavec},t)$ is locally bounded with respect to
$\bfx,\hat\thetavec$ and uniformly bounded in $t$. Therefore the
difference $f(\bfx,\thetavec,t)-f(\bfx,\hat\thetavec,t)$ is
bounded. Furthermore, according to (\ref{varphi}), function
$\varphi(\psi)\in C^0$ and therefore, given that $\psi$ is
bounded, this function is bounded as well. Hence $\dpsi$ is
bounded and by applying Barbalat's lemma one can show that
$\psi(\bfx(t),t)\rightarrow 0$ at $t\rightarrow\infty$.

To compete the proof of the theorem  (property P5) consider the
difference $f(\bfx,\thetavec,t)-f(\bfx,\hat\thetavec,t)$. Let
function $\varphi\in C^1$, function $f(\bfx,\thetavec,t)$ is
differentiable in $\bfx$, $\thetavec$; derivative $ \pd
{f(\bfx,\thetavec,t)}/{\pd t}$ is bounded uniformly in $t$;
function $\alphavec(\bfx,t)$ is locally bounded with respect to
$\bfx$ and uniformly bounded with respect to $t$, then $d/dt
(f(\bfx,\thetavec,t)-f(\bfx,\hat\thetavec,t))$ is bounded. On the
other hand there exists the following limit
\[
\lim_{t\rightarrow\infty}\int_0^{t}
(f(\bfx,\thetavec,\tau)-f(\bfx,\hat\thetavec,\tau))^{2}=\int_0^{\infty}
(f(\bfx,\thetavec,\tau)-f(\bfx,\hat\thetavec,\tau))^{2} \leq
\frac{1}{D}V_{\hat{\thetavec}}(\hat{\thetavec}(0),\thetavec)
\]
as $\int_0^{t}
(f(\bfx,\thetavec,\tau)-f(\bfx,\hat\thetavec,\tau))^{2}$ is
non-decreasing and bounded from above. Hence by Barbalat's lemma
it follows that
$f(\bfx,\thetavec,\tau)-f(\bfx,\hat\thetavec,\tau)\rightarrow 0$
as $t\rightarrow\infty$. Notice also that
$\psi(\bfx(t),t)\rightarrow 0$ as $t\rightarrow \infty$. Then
$\dpsi\rightarrow 0$ as $t\rightarrow\infty$. {\it The theorem is
proven.}

{\it Proof of Theorem \ref{exp_stability_theorem}.} Consider the
following integral\footnote{We substitute the arguments of the
functions $\dpsi(\cdot)$ and $\psi(\cdot)$ with $t$. This means
that we consider them as functions of time.}
$\int_0^{t}(\dpsi(\tau)+\varphi(\psi(\tau))^2 d\tau$. It was shown
in Theorem \ref{stability_theorem} proof that
$\int_0^{t}(\dpsi(\tau)+\varphi(\psi(\tau))^2 d\tau\leq
\frac{1}{2D}\|\hat{\thetavec}(0)-\thetavec\|^{2}_{\Gamma^{-1}}$
along system (\ref{system1}), (\ref{control}),
(\ref{algorithm_dpsi}) solutions. Let us define
$\mu(t)=\dpsi(t)+\varphi(\psi(t))$, or
\begin{eqnarray}\label{prop_exp_conv_dpsi}
\dpsi=-\varphi(\psi)+\mu(t),
\end{eqnarray}
where $\int_{0}^{\infty}\mu^2(\tau)d\tau\leq
\frac{1}{2D}\|\hat{\thetavec}(0)-\thetavec\|^{2}_{\Gamma^{-1}}$.
According to the theorem conditions, $\varphi(\psi)=K\psi$, it is
possible to derive the solution of equation
(\ref{prop_exp_conv_dpsi}) as follows
$\psi(t)=\psi(0)e^{-Kt}+\int_0^{t}e^{-K(t-\tau)}\mu(\tau)d\tau$.
Hence
\begin{eqnarray}\label{prop_exp_conv_psi}
|\psi(t)|&\leq& |\psi(0)|e^{- Kt} + \sqrt{\left(\int_{0}^t
e^{-K(t-\tau)}\mu(\tau)d\tau\right)^2}\leq |\psi(0)|e^{- Kt} +
\sqrt{\int_{0}^t
e^{-2K(t-\tau)}d\tau\int_0^{t}\mu^{2}(\tau)d\tau}\nonumber \\
&\leq & |\psi(0)|e^{-
Kt}+\frac{1}{2}\sqrt{\frac{1}{KD}\|\hat{\thetavec}(0)-\thetavec\|^{2}_{\Gamma^{-1}}}.
\end{eqnarray}
Property P6) is thus proven. In order to prove property P7)
consider
\[
\dot{\hat{\thetavec}}=\Gamma(\dpsi+\varphi(\psi))\alphavec(\bfx,t)=\Gamma
(f(\bfx,\thetavec,t)-f(\bfx,\hat{\thetavec},t))\alphavec(\bfx,t).
\]
Function
\begin{eqnarray}
& &
D_1|\alphavec(\bfx,t)^{T}(\hat{\thetavec}-\thetavec)|\leq|f(\bfx,\thetavec,t)-f(\bfx,\hat{\thetavec},t))|\leq
D|\alphavec(\bfx,t)^{T}(\hat{\thetavec}-\thetavec)|\nonumber\\
& &
\alphavec(\bfx,t)^{T}(\hat{\thetavec}-\thetavec)(f(\bfx,\hat\thetavec,t)-f(\bfx,\thetavec,t))>0
\  \forall \ f(\bfx,\thetavec,t)\neq
f(\bfx,\hat{\thetavec},t).\nonumber
\end{eqnarray}
Therefore, there exists $D_1\leq\kappa(t)\leq D$ such that
\[
\dot{\hat{\thetavec}}=-\kappa(t)\Gamma
\alphavec(\bfx,t)^{T}(\hat{\thetavec}-\thetavec)\alphavec(\bfx,t)=-\kappa(t)\Gamma
\alphavec(\bfx,t)\alphavec(\bfx,t)^{T}
(\hat{\thetavec}-\thetavec).
\]
Hence
\begin{eqnarray}\label{prop_exp_conv_dtheta}
\hat{\thetavec}(t)-\thetavec=e^{-\Gamma
\int_{0}^{t}\kappa(\tau)\alphavec(\bfx(\tau),\tau)\alphavec(\bfx(\tau),\tau)^{T}d\tau
}(\hat{\thetavec}(0)-\thetavec)
\end{eqnarray}
Consider the integral
$\Gamma\int_{0}^{t}\kappa(\tau)\alphavec(\bfx(\tau),\tau)\alphavec(\bfx(\tau),\tau)^{T}d\tau$
for $t>L$
\[
\Gamma\int_{0}^{t}\kappa(\tau)\alphavec(\bfx(\tau),\tau)\alphavec(\bfx(\tau),\tau)^{T}d\tau\geq
\Gamma D_1
\int_{0}^{t}\alphavec(\bfx(\tau),\tau)\alphavec(\bfx(\tau),\tau)^{T}d\tau,
\]
where $\alphavec(\bfx(t),t)$ is persistently exciting. For any
$t>L$ there exists integer $n\geq 0$ such that $t=nL+r$, $r\in R,
0 \leq r < L$. Therefore
\[
\Gamma D_1
\int_{0}^{t}\alphavec(\bfx(\tau),\tau)\alphavec(\bfx(\tau),\tau)^{T}d\tau\geq
\Gamma D_1 n \delta I \geq \left(\frac{\Gamma D_1
\delta}{L}t-I\right).
\]
Then taking into account (\ref{prop_exp_conv_dtheta}) one can
write
\begin{eqnarray}\label{prop_exp_conv_norm_theta}
\|\hat{\thetavec}(t)-\thetavec\|\leq \|e^{\left(-\frac{\Gamma D_1
\delta}{L}t+I\right)}\|\|\hat{\thetavec}(0)-\thetavec\|,
\end{eqnarray}
i. e. $\hat\thetavec(t)$ converges to $\thetavec$ exponentially
fast. It means that there exist positive constants $\lambda>0$,
$\lambda\neq K$ and $D_{\hat{\thetavec}}>0$ such that
$\|\hat{\thetavec}(t)-\thetavec\|\leq e^{-\lambda t}
\|\hat{\thetavec}(0)-\thetavec\| D_{\hat{\thetavec}}$. It follows
from Theorem \ref{stability_theorem} that $\psi(\bfx(t),t)$ is
bounded. In addition due to Assumption \ref{assume:psi} we can
conclude that $\bfx$ is bounded as well. By the theorem
assumptions function $\alphavec(\bfx,t)$ is locally bounded with
respect to $\bfx$ and uniformly bounded in $t$. Therefore, there
exists $D_{\alphavec}>0$ such that
$|\alphavec(\bfx,t)^{T}(\hat\thetavec(t)-\thetavec)|\leq
D_{\alphavec}\|\hat{\thetavec}(t)-\thetavec\|$. Taking into
account that $f(\bfx,\thetavec,t)-f(\bfx,\hat\thetavec,t)=\mu(t)$
and $|f(\bfx,\thetavec,t)-f(\bfx,\hat\thetavec,t)|\leq D
|\alphavec(\bfx,t)^{T}({\hat\thetavec}(t)-{\thetavec})|$ we can
derive from (\ref{prop_exp_conv_dpsi}) the following
 estimate
\begin{eqnarray}
|\psi(t)|\leq |\psi(0)|e^{-K t} + \|\hat{\thetavec}(0)-\thetavec\|
D_{\hat{\thetavec}}D_{\alphavec} D
\int_0^{t}e^{-K(t-\tau)}e^{-\lambda\tau}d\tau  \leq |\psi(0)|e^{-K
t} + \frac{D_{\hat{\thetavec}}D_{\alphavec}D}{K-\lambda}
\|\hat{\thetavec}(0)-\thetavec\| e^{-\lambda t}
\end{eqnarray}
{\it The theorem is proven.}

{\it Corollary \ref{cor:intervals} proof.} In order to prove the
corollary, we notice first that function $\sigma_j$ is equal to
unit for the following segments of the system solutions:
$\|\bfx(\bfx_0,t_0,t)-\bfc_j\|<r_j, \ \|\bfx_0-\bfc_j\|\leq
\delta_j$. Let us consider two cases: 1)
$\|\bfx(\bfx_0,t_0,t)-\bfc_j\|<r_j$ for any $t>t_0$, and 2)  for
any $t_0$ and  $\|\bfx(\bfx_0,t_0,t_0)-\bfc_j\|<r_j$ there exist
$t_1>t_0$ such that $\|\bfx(\bfx_0,t_0,t)-\bfc_j\|=r_j$.

In the first case Theorem \ref{exp_stability_theorem} explicitly
applies and the corollary follows automatically. In the second
case, we can derive from Theorem \ref{exp_stability_theorem} that
$\psi(\bfx,t)$ is bounded for every $t\in [t_0,t_1]$. Furthermore,
according to the properties of function $\psi(\bfx,t)$, it is
bounded in $t$  for every  $\bfx: \ \|\bfx-\bfc_j\|\leq \delta_j$.
Let  us denote this bound by symbol $\Delta_{\psi}$. Therefore,
according to Theorem \ref{stability_theorem} we can derive the
following estimate of $|\psi|_{\infty}$ for $t\in [t_0,t_1]$
\[
|\psi(\bfx,t)|\leq \Lambda
\left(Q(\Delta_{\psi})+\|\hat{\thetavec}(t_0)-\thetavec\|^2_{(4D\Gamma)^{-1}}\right)
\]
Given that norm
$\|\hat{\thetavec}(t_0)-\thetavec\|^2_{(4D\Gamma)^{-1}}$ is not
increasing, we can bound function $\psi(\bfx,t)$ for any time
moments $t: \sigma_j(t)=1$ as follows:
\[
|\psi(\bfx,t)|\leq \Lambda
\left(Q(\Delta_{\psi})+\|\hat{\thetavec}(0)-\thetavec\|^2_{(4D\Gamma)^{-1}}\right)
\]
On the other hand, due to the smoothness of function
$\psi(\bfx,t)$ and Assumption \ref{assume:steer_exists} one can
show that $\psi(\bfx,t)$ is bounded for every $t: \
\sigma_j(t)=0$. Hence, as follows from Assumption
\ref{assume:psi}, state $\bfx$ of the system is bounded. In order
to complete the proof we must show that
$\hat{\thetavec}(t)\rightarrow\thetavec$ as $t\rightarrow\infty$.
We have just shown that state $\bfx(t)$ is bounded. Then it is
bounded for those time intervals when $\sigma_j=1$ (i.e., when the
estimator is turned on). This implies that for any
$k=1,2,...,\infty$ the difference $t_k'-t_k>\delta_t>0$ (i. e.,
the time when the estimator is turned on is bounded from below).
Therefore, assuming that $L$ is sufficiently small (for instance,
$L<\delta_t/2$) and applying the same arguments as in the proof of
Theorem \ref{exp_stability_theorem}, we can show that
\[
\|\hat{\thetavec}(t_k')-\thetavec\|\leq \|e^{\left(-\frac{\Gamma
D_1 \delta}{L}(n-1)\right)}\|\|\hat{\thetavec}(t_k)-\thetavec\|,
\]
where $t_k'=t_k+nL+r$, $0\leq r < L$. {\it The corollary is
proven.}

\bibliographystyle{plain}
\bibliography{identification_fifo_embed}

\begin{thebibliography}{10}

\bibitem{Armstrong_1993}
B.~Armstrong-Helouvry.
\newblock Stick silp and control in low-speed motion.
\newblock {\em IEEE Trans. on Automatic Control}, 38(10):1483--1496, 1993.

\bibitem{Bai_2003}
Er-Wei Bai.
\newblock Frequency domain identification of hammerstein models.
\newblock {\em IEEE Trans. on Automatic Control}, 48(4):530--542, 2003.

\bibitem{Bastin92}
G.~Bastin, R.R. Bitmead, G.~Campion, and M.~Gevers.
\newblock Identification of linearly overparametrized nonlinear systems.
\newblock {\em IEEE Trans. on Automatic Control}, 37(7):1073--1078, 1992.

\bibitem{Boskovic_1995}
J.D. Boskovic.
\newblock Stable adaptive control of a class of first-order nonlinearly
  parameterized plants.
\newblock {\em IEEE Trans. on Automatic Control}, 40(2):347--350, 1995.

\bibitem{Box}
M.~J. Box, D.~Davies, and W.H. Swann.
\newblock {\em Non-linear Optimization Techniques}.
\newblock Oliver and Boyd, 1969.

\bibitem{Canudas_1999}
C~Canudas~de Wit and P.~Tsiotras.
\newblock Dynamic tire models for vehicle traction control.
\newblock In {\em Proceedings of the 38th IEEE Control and Decision
  Conference}. 1999.

\bibitem{Cao_2003}
C.~Cao, A.M. Annaswamy, and A.~Kojic.
\newblock Parameter convergence in nonlinearly parametrized systems.
\newblock {\em IEEE Trans. on Automatic Control}, 48(3):397--411, 2003.

\bibitem{Hansel_1992}
H.~Sompolinsky D.~Hansel.
\newblock Synchronization and computation in a chaotic neural network.
\newblock {\em Phys. Rev. Lett.}, 68:718--721, 1992.

\bibitem{Abbott_2001}
P.~Dayan and L.F. Abbott.
\newblock {\em Theoretical Neuroscience: Computational and Mathematical
  Modeling of Neural Systems}.
\newblock MIT Press, 2001.

\bibitem{Enqvist_CDC2002}
M.~Enqvist and L.~Ljung.
\newblock Estimating nonlinear systems in a neighborhood of
  {L}{T}{I}-approximants.
\newblock In {\em In Proc. of the 41st IEEE Conference on Decision and
  Control}, pages 1005--1010. 2002.

\bibitem{Eykhoff}
P.~Eykhoff.
\newblock {\em System Identification. Parameter and State Estimation}.
\newblock Univ. of Techn. Eindhoven, 1975.

\bibitem{Garulli_2002}
A.~Garulli, L.~Giarre, and G.~Zappa.
\newblock Identification of approximated hammerstein models in a worst-case
  setting.
\newblock {\em IEEE Trans. on Automatic Control}, 47(12):2046--2050, 2002.

\bibitem{Hansen92}
E.~Hansen.
\newblock {\em Global Optimization Uzing Interval Analysis}.
\newblock Marcel Dekker, 1992.

\bibitem{Hindmarsh_and_Rose}
J.L. Hindmarsh and R.M. Rose.
\newblock A model of neuronal bursting using 3 coupled 1st order
  differential-equations.
\newblock {\em Proc. R. Soc. Lond.}, B 221(1222):87--102, 1984.

\bibitem{Johansen_95}
T.A. Johansen and B.A. Foss.
\newblock Identification of non-linear system structure and parameters using
  regime decomposition.
\newblock {\em Automatica}, 31(2):321--326, 1995.

\bibitem{Kirkpatrick83}
S.~Kirkpatrick, C.~Gelatt, and M.~P. Vecchi.
\newblock Optimization by simulated annealing.
\newblock {\em Science}, 220:671--680, 1983.

\bibitem{Ljung_99}
L.~Ljung.
\newblock {\em System Identification: Theory for the User}.
\newblock Prentice-Hall, 1999.

\bibitem{Morgan_77}
A.~P. Morgan and K.~S. Narendra.
\newblock On the stability of nonautonomous differential equations
  $\dot{\bfx}=[\mathbf{A}+\mathbf{B}(t)]\bfx$ with skew symmetric matrix
  $\mathbf{B}(t)$.
\newblock {\em SIAM {J}. {Control} and {O}ptimization}, 37(9):1343--1354, 1992.

\bibitem{Narendra89}
K.~S. Narendra and A.~M. Annaswamy.
\newblock {\em Stable Adaptive systems}.
\newblock Prentice--Hall, 1989.

\bibitem{Narendra_66}
K.~S. Narendra and P.G. Gallman.
\newblock An interative method for the indentification of nonlinear systems
  using a hammerstain model.
\newblock {\em IEEE Trans. on Automatic Control}, AC-11(7):546--550, 1966.

\bibitem{Pacejka91}
H.B. Pacejka and E.~Bakker.
\newblock The magic formula tyre model.
\newblock In {\em Proceedings of 1-st Tyre Colloquium, Delft, October 1991},
  pages 1--18. 1993.
\newblock Supplement to Vehicle System Dynamics, vol. 21.

\bibitem{Panteley_2001}
E.~Panteley, A.~Loria, and A.~Teel.
\newblock Relaxed persistency of excitation for uniform asymptotic stability.
\newblock {\em IEEE Trans. on Automatic Control}, 46(12):1874--1886, 2001.

\bibitem{Pawlak_91}
M.~Pawlak.
\newblock On the series expansion approach to the identification of
  hammerstatin systems.
\newblock {\em IEEE Trans. on Automatic Control}, 36(6):763--767, 1991.

\bibitem{I_Petersen_2003}
I.~Petersen, T.~Johansen, J.~Kalkkuhl, and J.~Ludemann.
\newblock Wheel slip control using gain-scheduled {LQ} - {LPV/LMI} analysis and
  experimental results.
\newblock In {\em Proceedings of {I}{E}{E} European Control Conference,
  Cambridge, UK, September 1--4}. 2003.

\bibitem{t_fin_formsA&T2003}
I.~Y. Tyukin.
\newblock Algorithms in finite form for nonlinear dynamic objects.
\newblock {\em Automation and Remote Control}, 64(6):951--974, 2003.

\bibitem{ECC_2003}
I.~Yu. Tyukin, D.~V. Prokhorov, and Cees van Leeuwen.
\newblock Finite form realizations of adaptive control algorithms.
\newblock In {\em Proceedings of {I}{E}{E} European Control Conference,
  Cambridge, UK, September 1--4}. 2003.

\bibitem{ALCOSP_2004}
I.~Yu. Tyukin, D.V. Prokhorov, and C.~van Leeuwen.
\newblock Adaptive algorithms in finite form for nonconvex parameterized
  systems with low-triangular structure, august 30 -- september 1.
\newblock In {\em Proceedings of the 8-th {I}{F}{A}{C} Workshop on Adaptation
  and Learning in Control and Signal Processing (ALCOSP 2004)}. Yokohama,
  Japan, 2004.

\bibitem{Verdult_2002}
V.~Verdult, L.~Ljung, and M.~Verhaegen.
\newblock Identication of composite local linear state space models using a
  projected gradient search.
\newblock {\em Int. Journal of Control}, 75(16/17):1385--398, 2002.

\bibitem{Wilde}
D.~J. Wilde and C.S. Beightler.
\newblock {\em Foundation of Optimization}.
\newblock Prentice-Hall, 1967.

\bibitem{Wu_1990}
C.~Wu, J.C. Houk, K.Y. Young, and L.E. Miller.
\newblock Nonlinear damping of limb motion.
\newblock In J.M. Winters and Woo S-L.Y., editors, {\em Multiple Muscle
  Systems}, pages 214--235. Springer-Verlag, 1990.

\bibitem{Zhang_96}
Y.~Zhang, P.~Ioannou, and C.~Chien.
\newblock Parameter convergence of a new class of adaptive controllers.
\newblock {\em IEEE Trans. on Automatic Control}, 41(10):1489--1493, 1996.

\end{thebibliography}

\end{document}